\input amstex
\input epsf
\documentstyle{amsppt}
\mag=1200
\widestnumber\key{KenVer}
\TagsOnRight
\raggedbottom

\define \X{\frak X}
\define \Y{\frak Y}
\define \a{\alpha}
\define \e{\varepsilon}

\define \ga{\gamma}
\define \de{\delta}

\define \la{\lambda}

\define \({\left(}
\define \){\right)}
\define \[{\left[}
\define \]{\right]}
\define \ov{\overline}

\define \n{\goth n}
\define \wt{\widetilde}
\define \Ker{\operatorname{Ker}}
\define \Tr{\operatorname{Tr}}
\define \BZ{\Bbb Z}
\define \BN{\Bbb N}
\define \BQ{\Bbb Q}
\define \BR{\Bbb R}
\define \BT{\Bbb T}
\define \Orb{\text{Orb}}

\redefine \t{\bold t}
\redefine \phi{\varphi}

\def\smfrac#1/#2{\leavevmode\kern-.1em
    \raise.5ex\hbox{\the\scriptfont0 #1}\kern-.1em
    /\kern-.15em\lower.25ex\hbox{\the\scriptfont0 #2}}

\topmatter

\abstract
We study the arithmetic codings of hyperbolic
automorphisms of the 2-torus, i.e. the continuous mappings
acting from a certain symbolic space of sequences with a finite alphabet
endowed with an appropriate structure
of additive group onto the torus which preserve this structure and
turn the two-sided shift into a given automorphism of the torus.
This group is uniquely defined by an automorphism,
and such an arithmetic coding is a homomorphism
of that group onto $\BT^2$.
The necessary and sufficient
condition of the existence of a bijective arithmetic coding is obtained; it
is formulated in terms of a certain binary quadratic
form constructed by means of a given automorphism. Furthermore,
we describe all bijective arithmetic codings in terms
the Dirichlet group of the corresponding quadratic field. The minimum
of that quadratic form over the nonzero elements of  the lattice
coincides with  the minimal possible order of the kernel
of a homomorphism described above.
\endabstract

\title Bijective arithmetic codings of hyperbolic automorphisms
of the 2-torus, and binary quadratic forms
\endtitle
\author Nikita SIDOROV and Anatoly VERSHIK
\endauthor
\rightheadtext{Bijective arithmetic codings of hyperbolic automorphisms}
\address Steklov Institute of Mathematics at St.~Petersburg, 27 Fontanka,
St.~Petersburg 191011, Russia
\endaddress
\email sidorov$\@$pdmi.ras.ru, vershik$\@$pdmi.ras.ru \endemail

\keywords Hyperbolic automorphism of the torus, binary quadratic form,
homoclinic point, bijective arithmetic coding, minimal arithmetic
coding
\endkeywords
\subjclass 28D05, 11E16
\endsubjclass

\endtopmatter

\document

\head 0. Introduction \endhead

In this work we continue studying the
symbolic dynamics of ergodic automorphisms of the 2-torus.
\footnote""{Supported by the RFBR grant 96-01-00676 and the INTAS-RFBR
grant 95-418}
The dynamics of automorphisms of the torus is related more
to number theory than to the general theory of dynamical systems. This is
why their coding should be considered as a number-theoretic problem.
This was the main idea of \cite{Ver2} and subsequent papers
(see \cite{Ver1} and references therein); recently
it was developed in \cite{KenVer} and later in the dissertation \cite{Leb}.
Recall that to the Markov coding of hyperbolic automorphisms of the torus
and more general hyperbolic dynamical systems a number
of classical works have been devoted, see, e.g., \cite{AdWe},
\cite{Sin}, \cite{Bow}, \cite{GuSi}.
These papers are accented on the structure of Markov partitions,
but without special interest to
the arithmetic structure. For more details and the history of the
problem see the recent survey \cite{Ad}.

The quadratic case is  studied in detail below, and one sees that the
relationship with the theory of quadratic extensions and binary integral
quadratic forms becomes even deeper than before. We set
certain natural requirements on a symbolic realization of a hyperbolic
automorphism of the 2-torus (more precisely, on the maximal commutative
subgroup of $GL(2,\BZ)$ containing this automorphism), see Problem~1
in Section~1. Furthermore, we give the necessary and sufficient condition
on the existence of a mapping from a symbolic compactum
onto the 2-torus which we call an {\it arithmetic coding}.
Namely, arithmetic coding is
a mapping acting from the symbolic compactum provided
with ``almost group" structure onto the torus as ``almost homomorphism"
of this structure to the torus as an additive group (see Section~1
for the precise definitions and axiomatics).
It is proved that each arithmetic coding is naturally parametrized by a
homoclinic point of a given automorphism.
In our considerations we use two-sided decompositions of the points
of $\BT^2$ whose one-sided restriction coincide with the well-known
{\it $\beta$-expansions} (see \cite{Pa}); however, the two-sided
version proves to lead to new effects and problems.

The symbolic compactum in question is either Markov, if the determinant of
the matrix specifying an automorphism equals $-1$, or sofic otherwise.
It is proved that in both cases the compactum, after a certain factorization
of sequences of zero measure, turns into a group in addition
(Proposition~1.4).

An arithmetic coding is a specific mapping from the fixed symbolic
compactum $\X$ onto the torus. This mapping can be considered as
{\it expansions of the points of $\BT^2$ into the two-sided convergent
series with respect to the orbit of an arbitrary homoclinic
point}. It has the following form:
$$
\phi_\t(\e)=\lim_{N\to+\infty}\(\(\sum_{-N}^N \e_n T^{-n}\t\)\mod\BZ^2\),
$$
where $\e$ is a sequence from the compactum $\X$ and
$\t$ is a homoclinic point for $T$ written in coordinates of $\BR^2$
(see Theorem~1.2 for more details). Such expansions initially appeared in
\cite{Ver2}, \cite{Ver3}.

We also give a criterion of
the existence of a {\it bijective} arithmetic coding (see Theorem
below). In the case, where for a given
automorphism there is no bijective arithmetic coding, we present
a precise description of some minimal finite-sheeted covering of the torus.
A close connection with number theory that we mentioned above
is corroborated by the type of existence condition.

For the automorphism $T$ given by a matrix
$M_T=\(\smallmatrix a&b\\c&d\endsmallmatrix\)$
we define a very important {\it quadratic form
associated with $T$} by the formula
$$
f_T(x,y):=bx^2-(a-d)xy-cy^2.
$$
Let $\la$ be an eigenvalue of $M$, and $D$ be its discriminant. We
recall that the {\it Dirichlet group} $\Cal U_D$
of the quadratic field $\BQ(\sqrt D)$
is, by definition, the group of its units ($=$ units of its maximal
order), and by the classical theorem of number theory,
in our case the Dirichlet group is $\{(x+y\sqrt D)/2\}$, where $(x,y)$
is a solution of the Pell equation
$$
x^2-Dy^2=\pm4\tag0.1
$$
(see, e.g., \cite{BorSh} and \cite{Lev, vol. II, Chap. 2}). It is easy
to deduce from the cited theorem that if
a matrix $M$ is {\it primitive}, i.e. there is no matrix $K\in GL(2,\BZ)$
such that $M=K^n,\ n\ge2$, then $\Cal U_D=\{\pm\la^n\mid n\in\BZ\}
\simeq \BZ\times(\BZ/2\BZ)$.
Now we are ready to quote the essential part of the main result, see
Theorems~2.5 and 2.6 which concern the existence and properties
of the bijective arithmetic codings of $T$. Item~IV is taken
from Theorem~A.7 (see Appendix).

\proclaim{Theorem}
\roster
\item "I." The ergodic automorphism $T$ given by a matrix
$M_T=\(\smallmatrix a&b\\c&d\endsmallmatrix\)$ admits bijective arithmetic
coding
if and only its matrix $M_T$ is algebraically conjugate in $GL(2,\BZ)$ to the
companion matrix $C_{r,\sigma} =\(\smallmatrix r&1\\-\sigma&0\endsmallmatrix\)$
with $r=\Tr M_T$ and $\sigma=\det M_T$.
\item "II." A matrix $M_T$ is algebraically conjugate
to the corresponding companion matrix if and only if the Diophantine equation
$$
f_T (x,y)=\pm1
$$
is solvable.
\item "III." There exists a natural one-to-one correspondence
between the set of bijective arithmetic codings of $T$
and the Dirichlet group of the quadratic field $\BQ(\la)$ where $\la$
is an eigenvalue of $M_T$.
\item "IV." The existence of a point $(x,y)\in\BZ^2$ such that
the linear span for the orbit of $M_T$ of $(x,y)$ is equal
to $\BZ^2$, is equivalent to the algebraic conjugacy of
$M_T$ and the companion matrix $C_{r,\sigma}$. In particular,
necessarily $f_T(y,-x)=\pm1$.
\endroster
\endproclaim

More generally, a {\it minimal} arithmetic coding, i.e. a coding with
the minimal number of preimages, is naturally given by a solution
of the Diophantine equation $f_T (x,y)=\pm m$ with the minimal possible
positive $m$ (Theorem~3.5).

The precise axiomatic conditions on a symboilc realization of an automorphism
of the torus are as follows: the corresponding mapping from the
symbolic set of all sequences of nonnegative integers onto the 2-torus
should be a continuous homomorphism of semigroups turning the shift
into a given automorphism.
{\it A priori} it is not even clear, why so rigid conditions can be
satisfied. However, the fact that they really can,
yields a purely arithmetic interpretation
of a coding, namely, as two-sided convergent power series in powers
of the eigenvalue with a specifically chosen collection of digits and
Markov or {\it sofic} restrictions to their succession. This is nothing but
a two-sided generalization of the so-called $\beta$-expansions but
with essential sharpenings connected with the requirement of continuity
($=$ convergence).

A good deal of what was said above, might be extended to the general
case of a hyperbolic automorphism of $\BT^n,\ n\ge3$
whose principal eigenvalue is a PV number. Let us emphasize that
in higher dimensions in general it is not enough to consider natural
numbers as coefficients in the symbolic compactum; moreover, in \cite{KenVer}
it was shown that these coefficients could be algebraic numbers.
The condition of bijectivity is unknown for those cases.

Note also that for constructing examples which corroborate
some sharp estimates, we will use the facts from the theory of
indefinite quadratic forms contained, e.g., in the monograph
\cite{Cas1}, see Appendix. The relationship of this kind
of problems of dynamical systems theory with the geometry of numbers
and the theory of algebraic numbers becomes very important.
This link might be used in both directions.

\head 1. Basic notions and the main problem \endhead

\subhead 1.1. Basic notions and constructions \endsubhead
Let $\BT^2$ denote the 2-torus considered as the factor
$\BR^2/\BZ^2$. Let $T$ be an arbitrary group automorphism of $\BT^2$
given by a matrix
$\(\smallmatrix a&b\\c&d\endsmallmatrix\)\in GL(2,\BZ)$ which we will
denote by $M_T$. Let $r$ denote the trace of $M_T$, $\sigma$ stand for its
determinant. Suppose $T$ is hyperbolic,
which in the two-dimensional case is equivalent to the fact that none of the
roots of 1 belongs to the spectrum of $M_T$, i.e.
\roster
\item if $\sigma=-1$, then $r\neq0$;
\item if $\sigma=+1$, then $|r|\ge3$.
\endroster
The characteristic polynomial
of $M_T$ is $x^2-rx+\sigma$, and its dicriminant is
$D=r^2-4\sigma$. The eigenvalues of $M_T$ are $\frac12(r\pm\sqrt D)$.

Suppose $r$ to be positive; below we will prove that for our purposes
the study of the case of a negative trace will be immediate,
namely, we will consider the matrix $-M$ and easily
reformulate all the claims for it, see the end of Section~3.
Let $\la=\frac{r+\sqrt D}2>1$, and
and let $\ov\la$ denote the algebraic conjugate of $\la$, i.e.
$\ov\la=\sigma\la^{-1}=r-\la$.

We wish to consider symbolic codings
as appropriate expansions of the points of the torus in the sense
of some generalized ``number system" with natural coefficients.
Note that for multidimensional hyperbolic automorphisms
the coefficients are not necessarily naturals, but always elements
of a certain algebraic field, see \cite{KenVer}.
As a {\it primary symbolic set} of coefficients
for further codings we choose $\wt\X$ defined as the set of all
two-sided sequences with the coefficients
$\{\e_n\}_{-\infty}^\infty\in\prod_{-\infty}^\infty\BZ_+$
such that the series $\sum_{n=1}^\infty\e_n\la^{-n}$ and
$\sum_{n=1}^\infty\e_{-n}\la^{-n}$ converge.
We endow $\wt\X$ with the natural (weak) topology and
with coordinate-wise addition. It is obvious that $\wt\X$ is a
semigroup.

We call a sequence {\it finite}, if it contains only a finite number of
nonzero coordinates.
Let $\tau$ denote the two-sided shift on $\wt\X$,
i.e. $\tau\{\e_n\}=\{\e'_n\}$, where $\e'_n=\e_{n+1}$. We set up the main
problem of arithmetic coding axiomatically.

\definition{Definition} A one-sided sequence $(\e_1,\e_2,\dots)$ is
said to be {\it lexicographically less} than a sequence
$(\e'_1,\e'_2,\dots)$, if $\e_n<\e'_n$ for the least $n\ge1$
such that $\e_n\neq\e'_n$. Notation:
$(\e_1,\e_2,\dots)\prec_{\text{lex}}(\e'_1,\e'_2,\dots)$.
\enddefinition

\definition{Problem 1 (description of arithmetic codings)}
To describe all continuous semigroup homomorphisms
$\phi:\wt\X\to\BT^2$ which turn the shift $\tau$ into $T:\, \phi\tau=T\phi$.
For a given $\phi$ to find a {\bf closed, shift-invariant}
subset $\X$ of $\wt\X$ such that:
\roster
\item it is {\bf total}, i.e. together with a sequence
$\{\e_n\}_{-\infty}^\infty$ it contains
each sequence $\{\e'_n\}_{-\infty}^\infty$
such that
$(\e'_{N+k},\e'_{N+k+1},\dots)\prec_{\text{lex}} (\e_N,\e_{N+1},\dots)$
for some fixed $N\in\BZ$ and any $k\ge0$;
\item $\phi$ restricted to $\X$ is {\bf surjective} and
{\bf one-to-one on the set of finite sequences}
of $\wt\X$ belonging to $\X$ (the {\it section}
over the finite sequences).
\endroster
\enddefinition

\definition{Definition} For a hyperbolic automorphism $T$ of the 2-torus a
pair $(\phi, \X)$ defined in Problem~1, will be called
an {\it arithmetic coding} of $T$.
\enddefinition
We will see that such a coding exists for all hyperbolic automorphisms,
the compactum $\X$ depending on the spectrum of $M_T$ (not on $\phi$).
So, sometimes by a coding we will imply a mapping $\phi$ only.
Furthermore, we will show that after small glueings $\X$ acquires the
structure of a group and in fact $\phi$ restricted to $\X$
is a group homomorphism. An arithmetic coding is not necessarily bijective
almost everywhere, moreover, sometimes there is no bijective arithmetic
coding for a given $T$ at all.

\definition{Problem 2 (bijective and minimal arithmetic codings)}
To give the necessary and sufficient condition of the existence of
a bijective a.e. (with respect to the Lebesgue measure on $\BT^2$)
arithmetic coding for a given automorphism $T$ and
to describe all bijective arithmetic codings (BAC). If a BAC does not
exist, to find a minimal arithmetic coding (MAC) defined as
a coding with the minimal possible number of preimages and
to describe all such codings.
\enddefinition
We are going to solve Problem~1 in this section and to devote
two subsequent ones to Problem~2.

Let $\X_r$ denote the stationary {\it Markov compactum} with the
state space $0,1,\dots,r$ and the pairwise restrictions
$\{\e_n=r\Rightarrow \e_{n+1}=0,\ n\in\BZ\}$,
and the {\it sofic compactum}
$\Y_r=\{\{\e_n\}_{-\infty}^\infty: 0\le\e_n\le r-1,\
(\e_n\dots\e_{n+s})\neq(r-1)(r-2)^{s-2}(r-1)$ for any
$n\in\BZ$ and any $s\ge2\}$. Each of these compacta is the
$\beta$-{\it compactum} for $\beta=\la$. Let us give the corresponding
definition (see \cite{Pa}).

\definition{Definition} Let $\beta>1$, and
$1=d_1\beta^{-1}+d_2\beta^{-2}+\dots$, where
$d_1(\beta)=[\beta], d_2(\beta)=[\beta\{\beta\}],\dots$. Then by definition,
$X_\beta=\{\{\e_n\}_{-\infty}^\infty:(\e_n,\e_{n+1},\dots)\prec_{\text{lex}}
(d_1,d_2,\dots),\allowmathbreak n\in\BZ\}$.
The set $X_\beta$ endowed with the weak
topology is called the {\it $\beta$-compactum}.
\enddefinition

We need to recall one more classical definition.

\definition{Definition} Let $T$ be a hyperbolic automorphism of the torus.
A point $x$ is called {\it homoclinic} (to zero), if
$T^n(x)\to\bold 0$ as $n\to\pm\infty$.
\enddefinition
The equivalent definition is that $x$ belongs to the intersection of
the leaves of the stable and unstable foliations for $T$
going through $\bold 0$.

A suitable way of obtaining all homoclinic points for a given
automorphism was proposed in \cite{Ver3}.
Let $T$ be a hyperbolic automorphism of $\BT^k$ (not necessarily
two-dimensional). Consider the linear subspace of $\BR^k$ containing
the leaf of the unstable foliation going through $\bold 0$.
Then the projection of a point of the lattice $\BZ^k$
to this subspace along the direction of the stable foliation
taken modulo $\BZ^k$ is always a homoclinic point for $T$, and
and any homoclinic point can be obtained in such a way
(see \cite{Ver3} for more detail). For the two-dimensional
case these considerations yield the following complete description
of the homoclinic points.

\proclaim{Lemma 1.1} For the hyperbolic automorphism $T$ given by a
matrix $M_T$ its any homoclinic point $\t$ is parametrized by a pair
$(u,v)\in\BZ^2$ as follows: $\t=(\xi,\eta)$, where
$$
(\xi,\eta)=\(\frac{v+n\la}{\sqrt D},\frac{u+k\la}{\sqrt D}\)\mod\BZ^2,\tag1.1
$$
and
$$
\(\matrix n \\ k\endmatrix\)=-\det M_T\cdot M_T
\(\matrix v \\ u\endmatrix\).\tag1.2
$$
\endproclaim
\demo{Proof} Let $M_T=\(\smallmatrix a&b\\c&d\endsmallmatrix\)$.
From the general approach decribed above it follows that
to obtain any homoclinic point, one needs to consider
the projection $(\xi_0,\eta_0)$ of a certain point $(n,k)\in\BZ^2$ onto the
eigenline $y=\frac{\la-a}b x$ along the eigenline $y=\frac{\ov\la-a}b x$;
then this homoclinic point is $(\xi_0,\eta_0)$ modulo $\BZ^2$.
Solving the equation $\frac{\xi_0-n}b=\frac{\eta_0-k}{\ov\la-a}$ together
with $b\eta_0=(\la-a)\xi_0$, we get
$$
\xi_0=\frac{-dn+bk+n\la}{\sqrt D},\,\,\eta_0=\frac{cn-ak+k\la}{\sqrt D}
$$
(in view of the relations $\la+\ov\la=r,\la-\ov\la=\sqrt D$). Setting
$v:=-dn+bk,u:=cn-ak$, we complete the proof, because
$\(\smallmatrix -d&b\\c&-a\endsmallmatrix\)=-\det M_T\cdot M_T^{-1}$.\qed
\enddemo
\remark{Remark {\rm1}} As we see, there is a natural one-to-one
correspondence between the homoclinic points of $T$ and the
projections of the integral points onto the eigenline of $M_T$ being the
leaf of the unstable foliation going through $\bold0$
along its another eigenline. This fact
gives us an occasion to use below coordinates in $\BR^2$ for
the homoclinic points of $T$, which looks more natural.
\endremark
\remark{Remark {\rm2}} The purpose of such a choice of parameters in
Lemma~1.1 will become clear below, see Theorem~3.1.
\endremark

In the two-dimensional case that we are dealing with,
this approach leads to the fact that the group of homoclinic points
for $T$ is isomorphic to $\BZ[\la]$, this is why it will be convenient
to treat norm in $\BQ(\la)$ by means of homoclinic points.

\definition{Definition} Let $\|x\|:=\min\,\{|x-n|:n\in\BZ\}$.
A two-sided series of reals $\sum_{-\infty}^\infty w_n$
is said to converge to $w\in[0,1)$ {\it modulo 1}, if
$\|\sum_{n=-k}^l w_n-w\|\to0$ as $k,l\to+\infty$. The convergence
of a pair of series {\it modulo} $\BZ^2$ to a point of the torus
means the convergence of each coordinate modulo 1.
Besides, we will use the following notation:
$\(\sum\limits_{-\infty}^\infty w_n\)\(\smallmatrix \xi\\ \eta
\endsmallmatrix\)\mod\BZ^2:=\lim\limits_{N\to+\infty}
\(\sum\limits_{-N}^N\xi w_n,\sum\limits_{-N}^N\eta w_n\)
\mod\BZ^2$. Besides, by mutliplication of a homoclinic
point by some integer we imply the operation of multiplication
in the planar coordinates with (if necessary) further return
to the toral coordinates.
\enddefinition

\proclaim{Theorem 1.2}
\roster
\item "(I)" The set of arithmetic codings $(\phi,\X)$
of a hyperbolic automorphism $T$
is in a one-to-one correspondence with the points homoclinic to zero.
This correspondence is given by the formula for $\phi=\phi_\t$:
$$
\aligned
\phi_\t(\e)=&\lim_{N\to+\infty}\(\(\sum_{-N}^N \e_n T^{-n}\t\)\mod\BZ^2\)\\
           =&\(\sum_{n=-\infty}^\infty\e_n\la^{-n}\)\t\mod \BZ^2,
\endaligned\tag1.3
$$
where $\t$ is a homoclinic point for $T$ (in
coordinates of $\BR^2$). Conversely, if $\t$ is a homoclinic
point for $T$, mapping~{\rm(1.3)} is a convergent series and specifies an
arithmetic coding for $T$. Furthermore, if $\t\neq\bold0$, then
$\phi_\t$ is surjective.
\item "(II)" For each coding, a shift-invariant subset $\X$ satisfying
the second condition of Problem~1,
i.e. the surjectivity of the mapping $\phi_\t|_{\X}$ and its bijectivity
on the set of finite
sequences belonging to $\X$, is the stationary Markov compactum $\X_r$
for $\det M_T=-1$ or the sofic compactum $\Y_r$ for $\det M_T=+1$.
\endroster
\endproclaim
\demo{Proof} (I) Let $\phi$ satisfy the conditions of Problem~1.
We denote $u_k=\tau^k(u_0)$, i.e the sequence having 1 at the $(-k)$'th place
and zeroes at all other places. We set $\t=(\xi,\eta):=\phi(u_0)$.
By virtue of the continuity of the mapping $\phi$ and the fact that
$u_k\to0,\ k\to\pm\infty$, we have $T^k\t\to\bold0,\ k\to\pm\infty$,
whence by definition, $\t$ must be a homoclinic point. Hence
$\phi(u_k)=\la^k\t\mod\BZ^2$ for all $k\in\BZ$.

Consider now an arbitrary finite sequence $\e=\{\e_n\}\in\wt\X$.
By the additivity of $\phi$, we have
$\phi(\e)=\sum_n\e_n\phi(u_{-n})\mod\BZ^2$, whence
$$
\phi(\e)=\(\sum_{n=-\infty}^\infty\e_n\la^{-n}\)
\pmatrix \xi \\ \eta \endpmatrix\mod \BZ^2.
$$
We can now extend the mapping $\phi$ by continuity to all sequences
$\e\in\wt\X$, because since $\t$ is a homoclinic point,
$\la^N\t\to\bold0$ as $N\to\pm\infty$ with exponential rate of convergence,
whence $(\sum_{|n|>N}\e_n\la^{-n})\t\mod\BZ^2\to\bold 0$ as $N\to+\infty$
for any sequence $\{\e_n\}_{-\infty}^\infty\in\wt\X$.
Thus, if a mapping $\phi$ is an arithmetic coding, it must have form~(1.3).

Conversely, let a mapping $\phi_\t$
from $\wt\X$ onto the 2-torus be specified by
formula~(1.3) with $\t=(\xi,\eta)$ being a homoclinic point written
in coordinates of $\BR^2$.
The convergence of the series involved follows from the definition
of $\wt\X$.
We need to check that $\phi_\t$ is additive, continuous and turns the shift
into $T$. The {\bf additivity} of $\phi_\t$ on $\wt\X$ is a consequence of
its obvious coordinate-wise additivity. To prove its {\bf continuity},
consider two sequences $\e$ and $\e'$ such that $\e_n=\e'_n$ for
$|n|\le N$. Then
$\phi_\t(\e')-\phi_\t(\e)=
\bigl(\sum_{|n|>N}(\e'_n-\e_n)\bigr)\t\to\bold0$ as $N\to+\infty$.
As $\phi_\t$ is continuous, it suffices to verify the {\bf relation}
$\phi_\t\tau=T\phi_\t$
on the set of finite sequences. Let $\e\in\wt\X$ be finite;
then $\phi_\t\tau(\e)=\la\t\mod\BZ^2$, and $T\phi_\t(\e)=T\t=\la\t\mod\BZ^2$,
because $\t$ is homoclinic. Finally, let $\t\neq\bold0$. To prove the
{\bf surjectivity} of the mapping $\phi_\t$, we rewrite formula~(1.3) in the
form
$$
\phi_\t(\e)=\lim_{N\to+\infty}\(\sum_{n=-N}^\infty\e_n\la^{-n}\)
\t\mod \BZ^2.
$$
Thus, the image $\phi_\t(\wt\X)$ is the closure of the leaf of the unstable
foliation going through $\bold0$, whence this image is $\BT^2$, because
the leaf has irrational slope and thus is dense.
\smallskip
(II) Suppose now $\phi_\t|_\X$ (we will keep the same notation
$\phi_\t$ for this restriction) to be bijective on the set of finite
sequences for some shift-invariant subset $\X$ of $\wt\X$.
Our goal consists in showing that $\X=\X_r$ in the case $\sigma=-1$
or $\Y_r$ otherwise. Let for simplicity $\sigma=-1$. We first
prove the inclusion $\X\subset\X_r$. Let, on the contrary, $\X$
contain a sequence $\e$ not lying in the Markov compactum; then
there exists $k\in\BZ$ such that either $\e_k=r,\e_{k+1}\ge1$
or $\e_k>r$. Recall that by our assumption, $\X$ is total, and
thus, the second case implies the first one. Therefore, $\X$
contains the sequence
$(\dots,0,0,\dots,0,r,1,0,0,\dots)$.
The existence of such a sequence contradicts the assumption that
$\phi_\t$ is bijective on the finite sequences, because
$ru_n+u_{n+1}=u_{n-1}$.

To prove the inverse inclusion, suppose $\X\subsetneqq\X_r$. By
the stationarity and closeness of $\X$, this means that there exists a
cylinder $\{\e_0=i_0,\dots,\e_n=i_n\}$ belonging to $\X_r\setminus\X$
together with all its shifts. Below it will be shown that there exists
an ergodic measure $\mu$ on $\X_r$ such that $\phi_\t(\mu)$ is the
two-dimensional Lebesgue measure. Hence by the ergodic theorem,
$\mu(\X)=0$, and by the fact that
any mapping $\phi_\t:\X_r\to\BT^2$ of the form~(1.3)
is bounded-to-one (see Proposition~1.4 below), the Lebesgue
measure of the image $\phi_\t(\X)$ would be equal to 0, which contradicts
the surjectivity of $\phi$.\qed
\enddemo

\subhead 1.2. Group interpretation of an arithmetic coding \endsubhead
The Markov (or sofic) compactum defined above does not form a subsemigroup
of the semigroup $\wt\X$. Neverthelss, we can introduce a group structure
after certain small glueings of some sequences.
Let $\X_r^{(0)}$ and $\Y_r^{(0)}$
denote the subsets of $\X_r$ and $\Y_r$ respectively consisting of all
finite sequences. We have shown above that $\X_r^{(0)}$ is in fact
the factor of the semigroup
$\X^{(0)}=\{\sum_{|n|<N}\e_nu_n\mid \{\e_n\}\in\wt\X,\ N\in\Bbb N\}$
with respect to the stationary recurrence relations
$\{u_{n-1}=ru_n+u_{n+1},\ n\in\BZ\}$. Similarly, $\Y_r^{(0)}$
is the factor of $\X^{(0)}$ with respect to the relations
$\{u_{n-1}=ru_n-u_{n+1},\ n\in\BZ\}$.

It is well-known that in both cases in question the finite sequences
themselves form an additive semigroup
(see, e.g., \cite{FrSa}). Our goal now consists in assigning
the structure of an additive group to the whole symbolic compacta.
To do this, we first give the well-known definition of normalization
(see \cite{Fr}).
\definition{Definition} Let $x\in\prod_1^{\infty}\BZ_+,
\ x=\{x_k\}_{k=1}^\infty$;
we define $c(x)=\sum_{k=1}^\infty x_k\la^{-k}=
\sum_{k=k_0}^\infty \e_k\la^{-k}$, where $\{\e_k\}$ is the $\beta$-expansion
of $c(x)$, i.e. the expansion whose digits are given by the greedy algorithm.
Thus, $(\e_{-k_0},\e_{-k_0+1},\dots)$ belongs to the symbolic compactum
$\X_r$ or $\Y_r$ respectively.
We define
$$
\n(x):=\{\e_k\}_{k=k_0}^\infty.
$$
The operation $\n$ is called the {\it normalization} of a sequence.
\enddefinition
With the help of normalization we can now define addition and subtraction
on the symbolic compacta. Let the elements of a sequence $x$ from the
definition of normalization are uniformly bounded, for instance,
$0\le x_k\le 2r$. Then for the cases in question (i.e. for the quadratic
units) it is
known that similarly to addition, the normalization of a finite sequence
is also finite and the carry to both sides is uniformly bounded, see
\cite{FrSo}.
Thus, it is easy to define the {\it two-sided normalization} of
almost every sequence with respect to any shift-invariant
measure $\mu$ being positive on each cylinder.
Namely, by the result of Frougny and Solomyak cited above, there
exists $L=L(\la)\in\BN$ such that the one-sided normalization of
any sequence with coefficients less than or equal to $2r$ which
has infinitely many blocks $(0\dots0)$ ($L$ times) is blockwise.
Thus, one can define the two-sided normalization for any sequence
containing this block infinitely many times to both sides from
the zero place. Note that the existence of such a block is not
necessary but sufficient. For more details see \cite{SidVer},
where the precise procedure was described in the case $\la=\frac{\sqrt5+1}2$.

The below theorem-definition is based on the following consideration.
We need to define subtraction on $\X$, specifically, the operation
$i:\e\mapsto-\e$. To do it, we are going to find for each of the
compacta invloved a sequence which is naturally identified with the
zero sequence in the sense of the arithmetic. For this goal we
consider different representations in $\X$ of the
elements $u_n=\tau^n(u_0)$ and easily see that for the Markov
compactum $\X_r$,
$$
u_n=ru_{n-1}+ru_{n-3}+\dots,
$$
and for the sofic compactum $\Y_r$,
$$
u_n=(r-1)u_{n-1}+(r-2)u_{n-2}+(r-2)u_{n-3}+\dots,
$$
whence for $\X_r$ the sequences $(\dots,r,0,r,0,\dots)$ are by continuity
identified with the zero sequence, the same is true for
the sofic case with the sequence $(\dots,r-2,r-2,r-2,\dots)$.
Our idea is to define the operation $i(\e)$ for the Markov compactum as
the normalization of the sequence defined as $\e_n'=r-\e_n$,
similarly, as $\e'_n=2(r-2)-\e_n$ for the sofic compactum, i.e.
to define $-\e$, we subtract $\e$ from the sequence whose
normalization is the zero sequence. Here is the precise claim.

\proclaim{Theorem-Definition} {\rm(concerning the group structure on $\X$).}
Let $\X$ denote one of the compacta $\X_r$ or $\Y_r$. We define the
operations of summation and turning to the inverse element in addition
in $\X$ as follows: let $\e$ and $\e'$ belong to the compactum
$\X$; the sequence $x=\{x_k\}_{-\infty}^\infty$ is defined
as $x_k=\e_k+\e'_k$. Then the sum of $\e$ and $\e'$ is by definition
the two-sided normalization of $x$.
To define $-\e$, consider the sequence $y=\{y_k\}$
with $y_k=r-\e_k$ for the Markov case and $y_k=2(r-2)-\e_k$
for the sofic case. By definition, $-\e$ is the two-sided normalization of
$y$. Both operations are well defined for a.e. sequence (or pair
of sequences for summation) with respect to any Borel measure
which is positive on each cylinder in $\X$ (respectively with respect
to the square of such a measure for addition).
\endproclaim
\demo{Proof} By the above, the sum of two sequences is well defined for
any pair $(\e,\e')$ such that both contain the block $(0\dots0)$
($L$ times) at
the same place infinitely many times to both sides from the zero coordinate.

The operation $i:\e\mapsto-\e$ in the sofic compactum $\Y_r$ is well defined
for $\e$ which has the block $(r-1,r-3,r-2,r-2,\dots,r-2,r-2,r-3,r-1)$
of length $L+2$ infinitely many
times to both sides. Indeed, the operation $\e_k\mapsto 2(r-2)-\e_k$
turns this block into the block $(r-3,r-1,r-2,r-2,\dots,r-2,r-2,r-1,r-3)$
whose normalization is $(r-2,0,0,\dots,0,0,r-2)$ with $L$ zeroes.
Since a.e. sequence $\{\e_k\}$
has such a block infinitely many times to both sides,
the normalization of $\{2(r-2)-\e_k\}$ is blockwise, and therefore
is well defined.

Finally, in the Markov case with $r\ge2$ (the case
$r=1$ was considered in \cite{SidVer}) the operation $i$ is well
defined, for instance, for the sequences having the cylinder
$\{\e_k=r-1,\e_{k+1}=r\}$ infinitely many times to both sides.
Indeed, the two-sided normalization acts by changing any
triple $(l,r,k)\mapsto (l+1,0,k-1)$ for $l\le r-1,k\ge1$, whence,
as is easy to see, for the sequence $\{\e'_n\}$ with $\e_n'=r-\e_n$, the
two-sided normalization is independent for the pieces
$(\dots,\e'_{k-1},\e'_k)$ and $(\e'_{k+1},\e'_{k+2},\dots)$. Thus,
we split a.e. sequence $\e$ into such pieces, so that
the normalization of $\e'$ is blockwise.\qed
\enddemo

Now we are going to make the above claim more precise. We describe
all identifications in $\X$ which turn it into a group in addition.

\proclaim{Proposition 1.3} Let $\X'_r$ and $\Y'_r$ denote the factor
sets $\X_r/\goth R_1$ and $\Y_r/\goth R_2$, where \newline
(1) $\goth R_1$ is the identification of the pairs of sequences
$(*k,r,0,r,\dots)\sim(*k+1,0,\allowmathbreak0,0,\dots)$, and
$(\dots0,r,0,r,0,k*)\sim(\dots r,0,r,0,r-1,k+1*)$,
where $*$ denotes one and the same arbitary admissible tail, and
$0\le k\le r-1$.\newline
(2) $\goth R_2$ is the identification of the pairs of sequences
$(*k,r-1,r-2,r-2,\dots)\sim(*k+1,0,0,0,\dots)$, and
$(\dots0,r-2,r-2,r-2,r-1,k*)\sim(\dots 0,0,0,0,k+1*),\ 0\le k\le r-2$.
Then the factor sets $\X'_r$ and $\Y'_r$ are groups in addition.
\endproclaim
\demo{Proof} The calculations based on the
relations $u_{n-1}=ru_n+u_{n+1},\ n\in\BZ$ for the Markov case and
$u_{n-1}+u_{N+1}=(r-1)u_n+(r-2)u_{n+1}+\dots+(r-2)u_{N-1}+(r-1)u_N,\
n\in\BZ,N\ge n$, lead exactly to the identifications mentioned
in the claim of the proposition. We omit technical computations.
For more details see \cite{Ver1}, \cite{Ver2} for the Markov
case with $r=1$ (more general cases are similar in techniques).\qed
\enddemo

\remark{Remark {\rm1}} It is easy to see that according to the rule
of glueing given in Proposition~1.3, there are some sequences
which are identified with two or three other ones. For instance,
by continuity the zero sequence is identified with
the sequences $(\dots,0,r,0,r.0,r,0,\dots)$ and
$(\dots,r,0,r,0.r,0,r,0,\dots)$ in the Markov compactum, and with
the sequence $(\dots,r-2,r-2,r-2,\dots)$ in the sofic compactum.
\endremark
\remark{Remark {\rm2}} The group $\X_r'$ (resp. $\Y'_r$) is a compact Abelian
group, hence it possesses the Haar measure $\mu$, which by definition
is Borel and positive on each cylinder, i.e. satisfies the conditions
of Theorem-Definition. The natural projection $\X_r\mapsto\X'_r$
(resp. $\Y_r\mapsto\Y'_r$) as a map of measure spaces is an isomorphism
$\pmod0$, which follows from the nature of identifications.
\endremark
\remark{Remark {\rm3}} Let, as above, $u_k$ denote the sequence
having all zeroes except one unity at the $k$'th place.
The operation $\e\mapsto\e+u_k$ in the sense of group structure
defined above, is the {\it two-sided} version of {\it adic transformation}
(see \cite{Ver2}). It turns out that the ordinary adic transformation
generates the action of $\BZ$ on the one-sided $\beta$-compactum,
while the case in question the addition of finite sequences generate
the action of $\BZ^2$.
\endremark

Below we will need the following claim.

\proclaim{Proposition 1.4} Any mapping $\phi$ from the definition of
arithmetic coding is well defined on the factor sets and is a group
homomorphism of the groups $\X'_r$ (resp. $\Y'_r$) and $\BT^2$.
Any arithmetic coding as a mapping from
$\X_r$ (resp. $\Y_r$) onto $\BT^2$ is always $K$-to-$1$ a.e. with respect
to the measure $\mu$ for some natural $K$.
\endproclaim
\demo{Proof} The factor map $\phi':\X'_r (\Y'_r)\to\BT^2$ is well defined,
because by the definition of identifications (see Proposition~1.3),
$\phi(\e)=\phi(\e')$, if $\e$ is identified with $\e'$.
Furthermore, by the nature of the arithmetic in $\X$, we have
$\phi'(\e\pm\e')=\phi'(\e)\pm\phi(\e')$. The second claim
follows from the theorem on the homomorphic image of a group, from
which $(\phi')^{-1}(x)=(\phi')^{-1}(\bold0)+\e$, where $\e$ is a sequence in
the preimage of $x\in\BT^2$. Thus, $\#(\phi')^{-1}(x)\equiv\text{const}$.
\qed
\enddemo
\remark{Remark} The precise value of the function
$K=K(u,v)$ will be computed in Section~3.
\endremark

\head 2. Bijective arithmetic codings of automorphisms and the
associated binary quadratic form \endhead

Below we will see that sometimes there are no bijective arithmetic codings
of a given automorphism; however, even if they do exist for a certain
homoclinic point $\t$, it can happen that for another homoclinic point
the mapping $\phi$ is not bijective a.e. Here is the simplest example.

\example{Example} Consider the {\it Fibonacci} automorphism $\Phi$
given by the
matrix $\(\smallmatrix 1&1\\1&0\endsmallmatrix\)$. The corresponding
Markov compactum is
$\X_\Phi=\{\{\e_n\}: \e_n\in\{0,1\},\ \e_n\e_{n+1}=0,\ n\in\BZ\}$,
and $\la=\frac{\sqrt5+1}2$.
By Theorem~1.2, an arithmetic coding of $\Phi$ is a mapping $\phi_\xi$
from $\X_\Phi$ onto $\BT^2$ of the form
$$
\phi_\xi(\{\e_n\})=\(\sum_{n=-\infty}^\infty \xi\e_n\la^{-n}\)
\pmatrix 1\\ \la^{-1} \endpmatrix \mod \BZ^2.
$$
Usually, the coefficients $\xi\e_n$ assume the values 0 and 1
(see, e.g., \cite{Ber}). However, this mapping (i.e. $\phi_1$) from
$\X_\Phi$ onto the torus proves to be not bijective, but actually 5-to-1
a.e. The kernel of the group homomorphism $\phi_1':\X'_\Phi\to\BT^2$
is the group $\Cal K=\{0^\infty,(1.000)^\infty,(0.100)^\infty,(0.010)^\infty,
(0.001)^\infty\}$, where point denotes the border between negative
and nonnegative coordinates of a sequence.
Thus, the preimage of a.e. point of the torus consists
of five sequences, the difference of any two of them being equal to
one of the sequences in $\Cal K$, and the compactum $\X_\Phi$ is
splitted into five parts $X_1\cup\dots\cup X_5$ such that $\phi_1|_{X_k}$
is bijective a.e. for $1\le k\le5$.

At the same time, as will be shown below, for the automorphism $\Phi$
given by the companion matrix
$M_\Phi=\(\smallmatrix 1&1\\1&0\endsmallmatrix\)$ a bijective arithmetic
coding does exist, and the proper choice of coefficients
is $\xi\e_n\in\left\{0,\frac1{\sqrt5}\right\}$, i.e the mapping
$\phi_{1/\sqrt5}$. In Figure~1 we depict the images of the sets
$\{X_k\}_1^5$ under the mapping $\phi_{1/\sqrt5}$. Each of these images
is the square with the side $\frac1{\sqrt5}$.

\bigskip
\epsfxsize=9cm
\centerline{\epsfbox{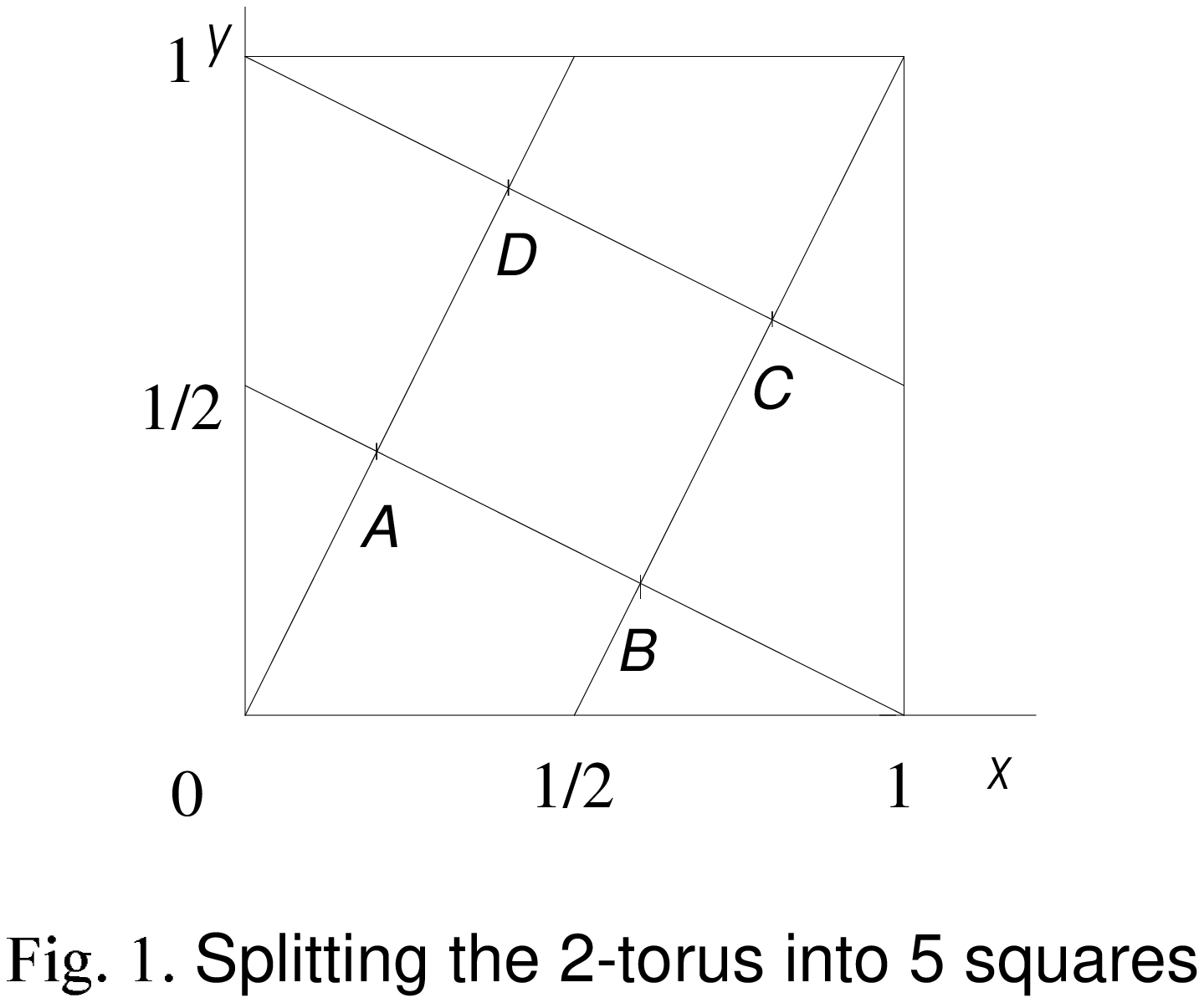}}
\bigskip\bigskip

The group
$\{O(0;0),A(\smfrac{1}/5;\smfrac{2}/5), B(\smfrac{3}/5;\smfrac{1}/5),
C(\smfrac{4}/5;\smfrac{3}/5), D(\smfrac{2}/5;\smfrac{4}/5)\}$
isomorphic to $\BZ/5\BZ$ is the image of the set $\Cal K$ under the mapping
$\phi_{1/\sqrt5}$. Note that the Fibonacci automorphism cyclically
moves the points of this group as follows: $\Phi:A\to B\to C\to D\to A$.

For a detailed study of the Fibonacci case and the proofs
see \cite{SidVer, sec.~1, item~1.6}.
\endexample

We will show that the condition on a homoclinic point $\t$
for the bijectivity a.e. of a mapping $\phi_\t$ given by formula~(1.3),
can be interpreted in terms of the area of
some fundamental domain. We begin with a class of matrices with
the simplest fundamental domain, namely with the case of {\it
companion} matrices. In this case $\t=(\xi,\pm\la^{-1}\xi)$ in
coordinates of $\BR^2$, and
the condition of bijectivity will be given in terms
of the algebraic norm of $\xi$.
Later it will be shown that the main result depends on the conjugacy
class in $GL(2,\BZ)$ and not on a matrix itself.

\subhead 2.1. Case of the companion matrix \endsubhead
We are going to
show that for the automorphism $T_{r,\sigma}$ given by
the {\it companion} matrix
$C_{r,\sigma} :=\(\smallmatrix r&1\\-\sigma&0\endsmallmatrix\)$
with $\sigma=\pm1$ and $r\in\Bbb N$ for $\sigma=-1$ and
$r\ge3$ for $\sigma=+1$, a BAC always exists
and that any such a coding is naturally parametrized by a unit
of the field $\BQ(\la)$. Note first that the vector
$\(\smallmatrix \la \\ -\sigma\endsmallmatrix\)$ is an eigenvector
of the matrix $C_{r,\sigma} $. Hence mapping~(1.3) in this case is given as follows:
$$
L_\xi(\{\e_n\})=\(\sum_{n=-\infty}^\infty \e_n\la^{-n}\)
\pmatrix \xi \\ -\sigma\xi\la^{-1} \endpmatrix \mod \BZ^2. \tag2.1
$$
To proceed, we need the precise description of possible values of $\xi$.
Recall that by the above, $(\xi,-\sigma\la^{-1}\xi)$ should be a homoclinic
point for $T_{r,\sigma}$. The following claim is a consequence of Lemma~1.1.

\proclaim{Lemma 2.1} The set of homoclinic points for the automorphism
$T_{r,\sigma}$ written in coordinates of $\BR^2$, is
$$
\left\{\(\frac{m+n\la}{\sqrt D},
-\sigma\la^{-1}\frac{m+n\la}{\sqrt D}\): (m,n)\in\BZ^2\right\}.
$$
\endproclaim
Thus, in formula~(2.1),
$$
\xi=\xi(m,n)=\frac{m+n\la}{\sqrt D}\tag2.2
$$
with $(m,n)\in\BZ^2$.\qed

Now our goal is to find among all $\xi$ of the form~(2.2) such that
$L_\xi$ is one-to-one a.e. We will see that actually these $\xi$ have
the minimal possible algebraic norm $N(\xi):=\xi\ov\xi$ in modulus, where
$\ov\xi$ is the algebraic conjugate of a quadratic irrational $\xi$.

\proclaim{Theorem 2.2} The automorphism
of the 2-torus $T_{r,\sigma}$ given by
the companion matrix
$C_{r,\sigma}:=\(\smallmatrix r&1\\-\sigma&0\endsmallmatrix\),\
\sigma=\pm1$,
admits a bijective arithmetic coding. If $\sigma=-1$ or $\sigma=+1,\,r\ge4$,
its BAC is always of the form
$$
L_\xi(\{\e_n\})=\(\sum_{n=-\infty}^\infty \e_n\la^{-n}\)
\pmatrix \xi \\ -\sigma\xi\la^{-1} \endpmatrix \mod \BZ^2,
$$
where $\xi=\frac{\pm\la^k}{\sqrt D}, k\in\BZ$.

The case $M=\(\smallmatrix 3&1\\-1&0\endsmallmatrix\)$ is specific.
Here $\la=\frac{3+\sqrt5}2$, and any BAC is of the form
$$
L_\xi(\{\e_n\})=\(\sum_{n=-\infty}^\infty \e_n\la^{-n}\)
\pmatrix \xi \\ -\xi\la^{-1} \endpmatrix \mod \BZ^2
$$
with $\xi=\frac{\pm\theta^k}{\sqrt D}, k\in\BZ,\
\theta=\sqrt\la=\frac{1+\sqrt5}2$.
\endproclaim
\demo{Proof} Let an arithmetic coding $L_\xi$ of $T_{r,\sigma}$
be written in the form~(2.1) with $\xi$ as in formula~(2.2).
Suppose first $\sigma=-1$.
Consider an arbitrary sequence $\{\e_n\}_{-\infty}^\infty\in\X$.
We split it into two pieces $\{\e_n\}_{-\infty}^0$ and $\{\e_n\}_1^\infty$
and define
$x_1(\{\e_n\}):=\sum_{k=1}^\infty \e_k\la^{-k},\ x_2=\sum_{k=0}^\infty
\e_{-k}(-\la)^{-k}$. It is a direct inspection that $x_1\in[0,1],\
x_2\in[-1,\la]$. Using the relation $\{\la^n\}=\{(-1)^{n+1}\la^{-n}\},\
n\ge0$, we make sure that
$\sum_{-\infty}^\infty\e_n\xi\la^{-n}=\xi x_1-\ov\xi x_2\mod1$
and similarly,
$\sum_{-\infty}^\infty\e_n\xi\la^{-n-1}=\xi\la^{-1}x_1+\ov\xi\la x_2\mod1$,
where, as above,
$\ov\xi$ denotes the algebraic conjugate of a quadratic irrational $\xi$.

Thus, we have the sequence of mappings
$$
\X_r @>F>> \Bbb R^2 @>b_\xi>> \Bbb R^2@>\pi>>\BT^2,
$$
where
$F(\{\e_n\})=(x_1,x_2)$, and
$b_\xi(x_1,x_2)=(\xi x_1-\ov\xi x_2,\,\xi\la^{-1}x_1+\ov\xi\la x_2)$,
i.e. the transfer to the eigenvector coordinates, and finally,
$\pi$ is the projection modulo the lattice $\BZ^2$. Thus, the mapping
$L_\xi$ is a factor map, i.e.
$$
L_\xi(\{\e_n\})=(\pi b_\xi F)(\{\e_n\}).
$$
By definition, the mapping $b_\xi F$ is always a bijection onto the image.
Note that since $(\e_0,\e_1)\neq(r,k)$ with $k\neq0$,
$F(\X_r)=
\Pi=([0,1]\times[-1,\la])\setminus([\la^{-1},1]\times[\la^{-1},\la])$,
i.e. the difference of rectangles (see Figure~2 below for the case
of the Fibonacci automorphism).
The area of $\Pi$ is $(\la+1)\la^{-1}+(1-\la^{-1})\la=\sqrt D$,
and the linear transformation
$b_\xi=\(\smallmatrix \xi & -\ov\xi \\ \la^{-1}\xi &
\la\ov\xi\endsmallmatrix\)$
from $\BR^2$ to $\BR^2$ has determinant $\sqrt D N(\xi)$,
where $N(\xi)=\xi\ov\xi$ is the algebraic norm of $\xi$.
Thus, the fundamental domain
$\Omega_\xi:=(b_\xi F)(\X_r)=b_\xi(\Pi)$ on the
plane has area $S=|DN(\xi)|$.

Recall that $\#L_\xi^{-1}(x)$ is one and the same for a.e. $x\in\BT^2$,
see Proposition~1.4. Thus, this capacity is necessarily equal to $S$,
and $L_\xi$ is a bijection a.e. if and only if the area $S$ of the
fundamental domain $\Omega_\xi$ equals 1, or equivalently, iff
$$
N(\xi)=\pm\frac1D. \tag2.3
$$

By Lemma~2.1, $\xi=\frac{m+n\la}{\sqrt D}$,
and the equation~(2.3) is equivalent to the Diophantine equation
$$
N(m+n\la)=\pm1.
$$
Therefore, as is well-known, $m+n\la$ is a unit of the ring $\BZ[\la]$
and thus, $m+n\la=\pm\la^k$ for some $k\in\BZ$ by virtue of the facts that
$\BZ[\la]$ is the maximal order of the field $\BQ(\la)$ and that
$\la$ is its main unit (see, e.g., \cite{BorSh}).
Let us recall that the above equation is in fact the condition on
a homoclinic point being the parameter of a coding.

The case $\sigma=+1$ is studied in the same way.
Since here $\ov\la=\la^{-1}$, we have
$x_1=\sum_1^\infty\e_k\la^{-k},\allowmathbreak
x_2=\sum_0^\infty\e_{-k}\la^{-k}$.
The set $\Pi$ here is the difference of the rectangles
$([0,1)\times[0,\la))\setminus((1-\la^{-1},1)\times(\la-1,\la))$.
The rest of the proof is the same, and we come to equation~(2.3).
Again, $m+n\la$ must be a unit of the ring $\BZ[\la]$, whence
$m+n\la=\pm\la^k$ if $r\ge4$, and $m+n\la=\pm\theta^k$ for $r=3$
with $\theta$ equal to the golden ratio.
\qed
\enddemo

Below we depict the fundamental domain $\Omega_\xi$ with $\xi=\frac1{\sqrt5}$
for the case of Fibonacci automorphism $\Phi$ (see Example above).

\bigskip
\epsfxsize=9cm
\centerline{\epsfbox{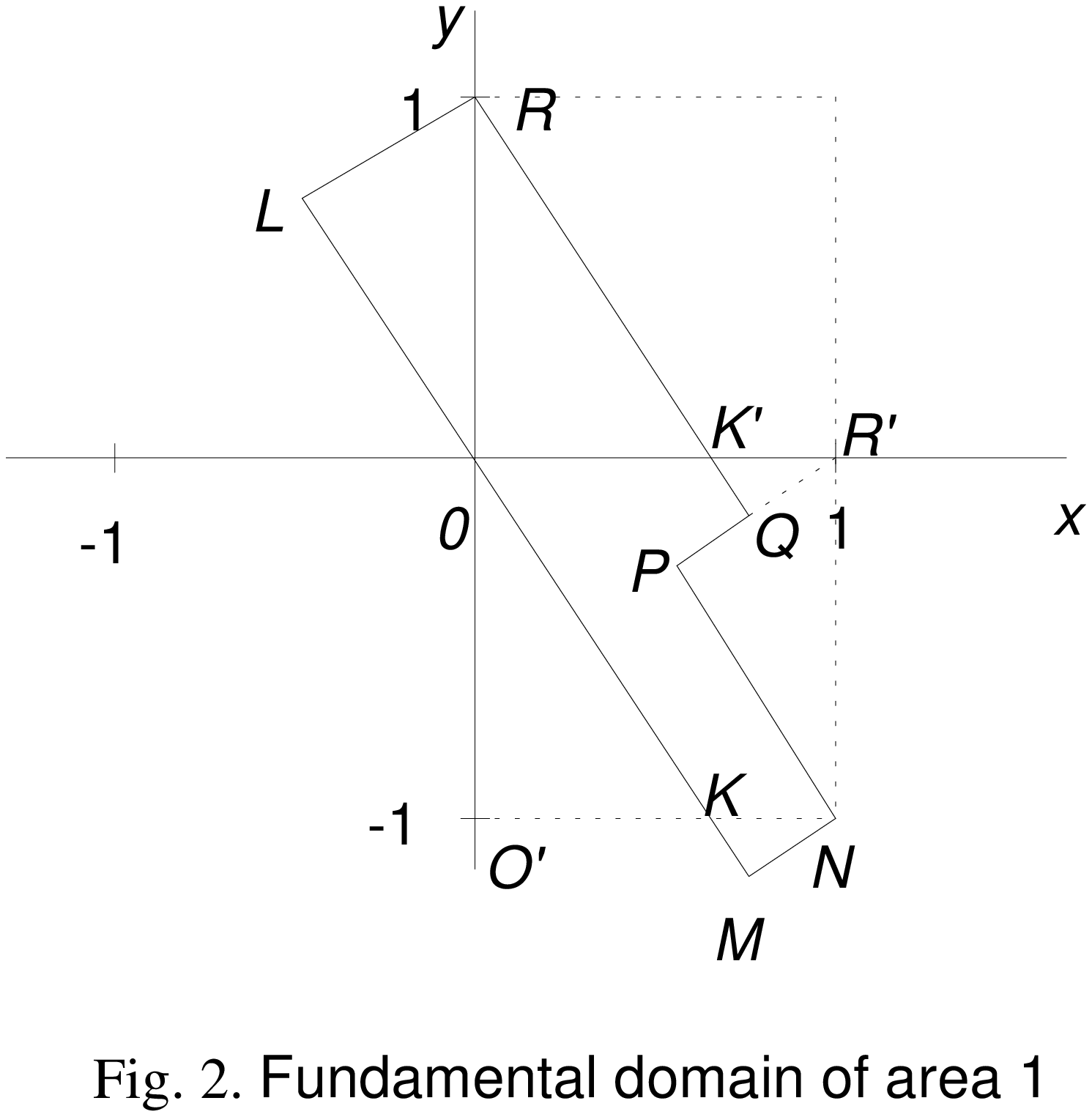}}
\bigskip\bigskip

It is visible from the figure that the fundamental domain is projected
modulo $\BZ^2$ onto the unit square. Indeed, consider the square
$O'OR'N=[0,1]\times[-1,0]$. The polygon $OKNPQK'$ lies inside the
square, and we project: triangle $ORK'$ onto $O'OK$, the triangle
$OLR$ onto $NPR'$, and finally, $MNK$ onto $QR'K'$.

\proclaim{Corollary 2.3} Any mapping $L_\xi$ from the symbolic compactum
onto the torus of the form~{\rm(2.1)} is $K$-to-$1$ with $K=|DN(\xi)|$.
\endproclaim

\remark{Remark} Let us give a geometric interpretation of the bijectivity.
We know that each parameter of an arithmetic coding of $T_{r,\sigma}$
is $\xi=\xi(m,n)=\frac{m+n\la}{\sqrt D}$, see formula~(2.2). Direct
computations show that those which yield a bijective arithmetic coding,
form the orbit (a kind of ``integral hyperbola")
$$
\left\{(m,n)\in\BZ^2:\(\matrix m\\n \endmatrix\)=
\pm C_{r,\sigma}^k\(\matrix 1\\0\endmatrix\)\,\,
\text{for some}\,\,k\in\BZ\right\},
$$
with the exception of the case $\sigma=+1,r=3$, when these integral
points form the orbit
$\left\{\pm \(\smallmatrix 1&1\\1&0\endsmallmatrix\)^k, k\in\BZ\right\}$.
\endremark
Recall that the homoclinic equivalence
relation on the torus is given as follows: two points $x$ and $y$
are said to be equivalent if $T^n(x-y)\to(0,0)$ as $|n|\to\infty$
(see, e.g., \cite{Gor}). Let us define the homoclinic equivalence
relation on the symbolic compactum.

\definition{Definition} A sequence in compactum $\X_r$
is called {\it homoclinic (to zero)} is its right and left tails are either
of the form $(0,0,0,\dots)$ or of the form $(r,0,r,0,\dots)$.
Similraly, a sequence in $\Y_r$ is called {\it homoclinic} if
its right and left tails are of the form $(0,0,0,\dots)$ or of the form
$(r-2,r-2,r-2,\dots)$.
Similarly to the ``toral" definition, we will say that two sequences
belong to the same {\it homoclinic class}, if their difference (which,
as we know, is well defined for a.e. pair, see Section~1) is a
sequence homoclinic to zero.
\enddefinition

\proclaim{Proposition 2.4} The image of the homoclinic class of a sequence
under a bijective arithmetic coding is the homoclinic class of
its image.
\endproclaim
\demo{Proof} By the above, after some identifications touching
sequences from one and the same homoclinic class, a BAC becomes
a complete bijection. Now the claim follows form the fact that a BAC
turns the sequences homoclinic to zero to the points homoclinic
to zero.\qed
\enddemo

\remark{Remark {\rm1}} For the case in question it seems more natural to
consider the following mapping which naturally generalizes the
one sided $\beta$-expansions to the two-sided ($=$ two-dimensional)
case:
$$
l(\{\e_n\})=\(\sum_{n=-\infty}^\infty \e_n\la^{-n}\)
\pmatrix -\sigma \\ \la^{-1} \endpmatrix \mod \BZ^2.
$$
It is a particular case of the more general mapping introduced and studied
in \cite{Ber}. Obviously, $l$ semiconjugates the shift $\tau$ and
the automorphism with the companion matrix $C_{r,\sigma} $; however,
from Corollary~2.3 it follows that the mapping $l$ is only $D$-to-1 a.e.
\endremark
\remark{Remark {\rm2}} For the case $\sigma=-1$ the mapping $L=L_{1/\sqrt D}$
in the form
$$
L(\e)=\sum_{k=1}^\infty \e_k\la^{-k}\cdot\pi_u(u_0)-
\sum_{k=0}^\infty \e_{-k}(-\la)^{-k}\cdot \pi_s(u_0) \mod \BZ^2
$$
with $u_0=\(\smallmatrix 0\\1\endsmallmatrix\)$ and $\pi_u,\,\pi_s$
being the projections on the leaves of the unstable and stable foliations
respectively,
was under consideration in the recent Ph.~D.~dissertation \cite{Leb}.
In particular, the author proved its bijectivity a.e. but did not
consider in detail its arithmetic properties.
\endremark

\subhead 2.2. General case \endsubhead
Return now to the general case of ergodic automorphism $T$ given by a matrix
$M_T=\(\smallmatrix a&b\\c&d\endsmallmatrix\)\in GL(2,\BZ)$. We begin with
two necessary definitions.

\definition{Definition} Two matrices $M_1$ and $M_2$ will be called
{\it algebraically conjugate},
if there exists a matrix $B\in GL(2,\BZ)$ such that
$BM_1B^{-1}=M_2$. We will write in this case $M_1\sim M_2$.
\enddefinition

\definition{Definition} The binary quadratic form
$f_T(x,y)=bx^2-(a-d)xy-cy^2$ will be called the {\it form associated
with an automorphism $T$}.
\enddefinition
\remark{Remark} Obviously, a binary integral quadratic form
$f(x,y)=\a x^2+\beta xy+\ga y^2$ is the form associated with some
automorphism if and only if $D(f)=r^2\pm4$ for some $r$, where
$D(f)=\beta^2-4\a\ga$ is the discriminant of the form $f$.
Since $D(f_T)=D$, we are dealing in fact with all forms with
the discriminant of the form $r^2\pm4>0$. The mapping
$\theta:T\mapsto f_T$ will be studied in detail in Appendix.
\endremark

\proclaim{Theorem 2.5}
\roster
\item "I." An ergodic automorphism $T$ admits bijective arithmetic coding
if and only its matrix $M_T$ is algebraically conjugate
to the
companion matrix $C_{r,\sigma}=\(\smallmatrix r&1\\-\sigma&0\endsmallmatrix\)$
with $r=\Tr M_T$ and $\sigma=\det M_T$.
\item "II." A matrix $M_T$ is algebraically conjugate 
to the corresponding companion matrix if and only if one of the equations
$$
f_T(x,y)=\pm1\tag2.4
$$
is solvable in $\BZ$.
Any matrix $B=\(\smallmatrix x&y\\z&t \endsmallmatrix\)\in GL(2,\BZ)$
such that $BM_TB^{-1}=C_{r,\sigma}$ is naturally parametrized by a solution of
Diophantine equation~{\rm(2.4)}, namely, $(x,y)$ is always a solution
of {\rm(2.4)}, and
$$
(z,t)=-\det M_T\cdot (x,y) M_T^{-1}.
$$
\endroster
\endproclaim
\demo{Proof} I. Suppose a matrix $B\in GL(2,\BZ)$
such that $BM_TB^{-1}=C_{r,\sigma}$
exists, and let $Q$ be the toral automorphism given by $B$. Let $L$
be a bijective a.e. mapping from relation~(2.1), say, for $\xi=1/\sqrt D$.
Recall that the compactum $\X$ is determined only by the spectrum of the
matrix specifying an automorphism, whence it is one and the same for $T$
and the automorphism given by $C_{r,\sigma}$. Consider the mapping
$\phi:=Q^{-1}L:\X\to\BT^2$. We have $\phi\tau=Q^{-1}L\tau=Q^{-1}T_{r,\sigma}=
TQ^{-1}L=T\phi$, and since $\phi$ is bijective a.e., it
is the desired BAC for $T$.

Conversely, let $T$ admit BAC, and $\phi$ be the corresponding
mapping from the symbolic compactum onto the torus. Consider
$Q:=L\phi^{-1}:\BT^2\to\BT^2$. It is well defined, because
if two sequences $\e,\e'$ belong to $L^{-1}(x)$ for some $x\in\BT^2$,
by the above, $\phi(\e)=\phi(\e')$. Thus, we make sure that by
definition of BAC, $Q$ is a group automorphism of the 2-torus,
hence, it is given by some matrix $B\in GL(2,\BZ)$. Since
$QTQ^{-1}=L\phi^{-1}T\phi L^{-1}=L\tau L^{-1}=T_{r,\sigma}$, we have
$C_{r,\sigma} =BM_TB^{-1}$.

II. It suffices to show that the solvability of
one of the Diophantine equations~(2.4)
is equivalent to the fact that $M_T\sim C_{r,\sigma}$.
Let $B=\(\smallmatrix x&y\\z&t\endsmallmatrix\)\in GL(2,\BZ)$ exist,
and $BM_T=C_{r,\sigma} B$. Suppose $\sigma=-1$.
We have thus the linear system
$$
\aligned
&z =-dx+cy\\
&t = bx-ay\\
&x = az+ct\\
&y = bz+dt,
\endaligned\tag2.5
$$
the last two equations being a consequence of the first two ones.
Hence this system together with the condition $\det B=\pm1$
yields the desired condition. For $\sigma=+1$ the first two equations
in formula~(2.5) are the same as for the previous case, so,
the argument is also the same.

Conversely, if the equation~(2.4) is solvable, then we take
some $x,y$ being its solutions and construct the matrix $B$ by
the equations for $z,t$ from formula~(2.5).
\qed
\enddemo

Recall that the Dirichlet theorem claims that given an automorphism $T$,
the group $\Cal D(T)$ defined as the set of all automorphisms of the
torus which  commute with $T$, has the form $\{\pm S^n,\ n\in\BZ\}$
for some primitive automorphism $S$. The following theorem shows that
in the Dirichlet group of $T$ only four primitive elements
$\pm S,\pm S^{-1}$ can admit BAC
(with the unique exclusion, when they are eight).

\proclaim{Theorem 2.6}
\roster
\item The automorphism $T$ admits BAC only if its matrix $M_T$
is primitive, with the exception of the case $M_T=K^2$, where
$K$ is algebraically conjugate to
$\(\smallmatrix 1&1\\1&0\endsmallmatrix\)$.
\item Let equation~{\rm(2.4)} be solvable for a given
automorphism $T$.
If $M_T$ is primitive, then there exists a homoclinic point
whose coordinates in $\BR^2$ are
$(\xi_0,\eta_0)$ such that any bijective arithmetic coding of $T$
is of the form
$$
\phi^\pm_k(\e)=\(\sum_{n=-\infty}^\infty \e_n\la^{-n}\)
\pmatrix \pm\xi_0\la^k \\ \pm\eta_0\la^k \endpmatrix \mod \BZ^2.\tag2.6
$$
If $M$ is not primitive, we have $\la=\frac{3+\sqrt5}2$, and
any BAC is of the form
$$
\phi^\pm_k(\e)=\(\sum_{n=-\infty}^\infty \e_n\la^{-n}\)
\pmatrix \pm\xi_0\theta^k \\ \pm\eta_0\theta^k \endpmatrix \mod \BZ^2,\tag2.7
$$
where $\theta=\sqrt\la=\frac{1+\sqrt5}2$.
\endroster
\endproclaim
\demo{Proof} (1) It is easy to compute that
$f_{T^2}=rf_T,\ f_{T^3}=(r^2+1)f_T$
and, more generally, $f_{T^n}=q_n(r)f_T$, where $q_n$ is a
polynomial of degree $n$ with nonnegative coefficients, odd
for $n$ odd and even for $n$ even, namely,
$q_n(r)=\frac1{\sqrt D}(\la^n-\ov\la^n)$. So, the form $f_{T^n}$ is
not {\it primitive} unless $n=2, r=1$, i.e. its coefficients are
not relatively prime, hence Diophantine equations~(2.4) have no
solutions. Thus, the unique companion matrix which is not primitive, is
$C_{3,1}\sim\(\smallmatrix 1&1\\1&0\endsmallmatrix\)^2$.
Now if $T$ admits bijective arithmetic coding, then by
Theorem~2.5, $M_T\sim C_{r,\sigma}$, and the first claim of the
theorem follows from the fact that the primitivity is an invariant
of algebraic conjugacy.

(2) Since equation~(2.4) is solvable in $\BZ$, there exists an infinite
number of different BAC's for $T$. Fix the notation
$\phi_0$ for one of them; let $\phi$ be
be an arbitrary BAC for $T$. Consider the mapping
$A:=\phi\phi_0^{-1}:\BT^2\to\BT^2$. It is well defined by the same arguments
as in the proof of Theorem~2.5. Obviously, $A$ is an automorphism of the
2-torus, and $A$ commutes with $T$. By the Dirichlet theorem cited above
and the primitivity of $T$, we have $\phi=\pm T^k\phi_0$ for some
$k\in\BZ$, whence if $\phi_0$ in formula~(1.3) is given by a homoclinic point
$(\xi_0,\eta_0)$, the mapping $\phi$ is given by
$(\pm\la^k\xi_0,\pm\la^k\eta_0)$.

Conversely, if a bijection a.e. $\phi$ is given by formula~(1.3) with some
$(\xi,\eta)$, the mapping $\phi'$ defined by the same formula with
$(\pm\la^{\pm1}\xi,\pm\la^{\pm1}\eta)$ is also a bijection, as
$\phi(-\e)=-\phi(\e),\ \phi(\tau^{\pm1}\e)=T^{\pm1}\phi(\e)$.
The argument for the exclusive case is the same with the
exception that here $\sqrt M$ is also a matrix in $GL(2,\BZ)$ and also
commutes with $M$.\qed
\enddemo
\remark{Remark} In Appendix we will give a simple example of a matrix
which is not conjugate to the corresponding companion matrix,
see ``Counterexamples".
\endremark
Thus, the bijective arithmetic codings in fact are naturally parametrized
by elements of the Dirichlet group of the field $\BQ(\sqrt D)$.

We finish the section by giving simple algebraic criteria for the
existence of a bijective arithmetic coding of a given automorphism
of the 2-torus.

\proclaim{Corollary 2.7} If two ergodic automorphisms $T_1$ and $T_2$
whose matrices have one and the same trace and discriminant,
both admit bijective arithmetic coding, then their matrices
are algebraically conjugate. Conversely, if $T_1$ admits BAC and
$M_{T_1}\sim M_{T_2}$, then so does $T_2$.
\endproclaim
\demo{Proof} It suffices to recall that both matrices should be
algebraically conjugate to the corresponding companion matrix
which is one and the same for both ones.\qed
\enddemo

\proclaim{Corollary 2.8} If $|b|=1$ or $|c|=1$, an automorphism $T$
with the matrix $\(\smallmatrix a&b\\c&d\endsmallmatrix\)$ admits
BAC.
\endproclaim
\demo{Proof} One of the equations~(2.4) has the trivial solution
$x=1,y=0$ or $x=0,y=1$ if $b=\pm1$ or $c=\pm1$ respectively.\qed
\enddemo

\head 3. Minimal arithmetic codings \endhead

We have already seen that sometimes an ergodic automorphism of the
2-torus does not admit BAC, so, it is meaningful to deal with
the notion of minimal arithmetic coding (MAC) introduced in Section~1.
Recall that a minimal arithmetic coding of an automorphism is,
by definition, a coding $\phi$ having the minimal possible number of
preimages.

Recall that by formulas~(1.1) and (1.2), any arithmetic coding
is parametrized by a pair $(u,v)\in\BZ^2$, and from Proposition~1.4
it follows that for any coding $(\phi,\X)$ a mapping $\phi:\X\to\BT^2$ is
$K$-to-1 a.e. The following theorem answers the question on the form of the
function $K=K(u,v)$.

\proclaim{Theorem 3.1} Let $T$ be the hyperbolic automorphism of the 2-torus
given by a matrix $M_T$. Then any arithmetic coding $\phi_\t$
of $T$ of the form~{\rm(1.3)}
with $\t=(\xi,\eta)$ being a homoclinic point defined by
formulas~{\rm(1.1)} and {\rm(1.2)}, is $K$-to-$1$ a.e. with
$$
K=K(u,v)=|f_T(u,v)|.
$$
\endproclaim
\demo{Proof}
Using the same arguments as in Theorem~2.2, we make sure that $K$ equals
the area of the fundamental domain and that this domain has area
given by the formula
$$
S=\sqrt D \left|\det \pmatrix \xi & -\ov\xi \\
\eta & -\ov\eta \endpmatrix\right|.\tag3.1
$$
Furthermore, from direct computations in formula~(3.1) which we omit
(in view of relation~(1.2)), it follows that
$$
K=S=|f_T(u,v)|.\tag3.2
$$
Thus, we proved the following theorem which describes explicitly,
in what way an arbitrary arithmetic coding is parametrized by
a homoclinic point.\qed
\enddemo

Let $m(T)$ denote the minimal possible number of preimages for an
arithmetic coding of $T$.

\proclaim{Corollary 3.2} The quantity $m(T)$ equals the integral minimum of
the associated form $f_T $. Any minimal arithmetic coding of a given
automorphism $T$ is naturally parametrized by a solution of the equation
$$
f_T(u,v)=\pm m, \tag3.3
$$
where $m=m(T)$.
\endproclaim

We are ready now to describe all possible minimal arithmetic codings
for a given automorphism more explicitly.

Let below $M'$ denote the transpose of $M$, and
$f_T$ stand also for the symmetric matrix of this quadratic
form, i.e
$$
f_T=\(\matrix b&\frac12(d-a)\\\frac12(d-a)&-c\endmatrix\).
$$

\proclaim{Lemma 3.3} We have
$$
M_Tf_TM_T' = \det M_T\cdot f_T,
$$
i.e. the change of variables given by the matrix
$M_T'$ turns the form $f_T$ into itself if $\det M_T=+1$ and into
$-f_T$ otherwise.
\endproclaim
\demo{Proof} Let, as above, $\sigma=\det M_T$. Then
$$
\align
M_Tf_TM_T'&=\pmatrix a&b\\c&d\endpmatrix
\pmatrix b&\frac12(d-a)\\\frac12(d-a)&-c\endpmatrix
\pmatrix a&c\\b&d\endpmatrix\\
&=\pmatrix \frac{br}2&-\frac{ar}2+\sigma\\\frac{dr}2-\sigma&-\frac{cr}2
\endpmatrix\pmatrix a&c\\b&d\endpmatrix=
\pmatrix \sigma b&\frac12\sigma(d-a)\\\frac12\sigma(d-a)&-\sigma c
\endpmatrix\\
&=\det M_T\cdot f_T.\qed
\endalign
$$
\enddemo

\definition{Definition} An integral change of variables which leaves
a binary integral quadratic form unchanged is called its {\it automorph}.
\enddefinition

Thus, if $\det M_T=+1$, then the transformation $M'_T$ is an automorph
of the form $f_T$.

Suppose from here on $M_T$ to be primitive. The following proposition
answers the question about the structure of the set of solutions
of equation~(3.3).

\proclaim{Proposition 3.4} Let $m$ denote the integral minimum
of the form $f_T$. The solutions of the equation~{\rm(3.3)}
are described as follows. The congruence
$$
n^2\equiv D\pmod{4m},
$$
is always solvable, and let $n$ be its minimum root, i.e.
$0\le n<2m$, and $l:=\frac{n^2-D}{4m}$.
Let $s$ stand for the number of distinct forms $[m,n,l]$ equivalent
to $f_T$.
Then there exists a finite collection of solutions of equation~{\rm(3.3)}
$(x^{(j)},y^{(j)}),\ 1\le j\le s$ such that any solution $(x,y)$ of
{\rm(3.3)} is of the form $(x,y)=\pm (x^{(j)}, y^{(j)})\cdot M_T^n$ for some
$n\in\BZ$ and $1\le j\le s$. Furthermore,
$(x^{(j)},y^{(j)})\neq\pm(x^{(i)},y^{(i)})\cdot M_T^n$ for $i\neq j$
and any integer $n$.
\endproclaim
\demo{Proof} We use the classical result on the structure of solutions
of a quadratic Diophantine equation (see \cite{Lev, vol.~II, Theorem~1-12},
by which if $(x,y)$ is a solution of the equation~(3.3), say, with $+m$,
then it leads to the series of solutions $\{V(x,y)'\}$, where $V$ is
an automorph of $f_T$. Besides, any solution of (3.3) is given
by such a series with a finite number of basis solutions. This
number is given exactly as in the claim. Furthermore,
dealing with $\pm m$, we see that
for our purposes we need to consider also the {\it anti-automorphs},
i.e. the transformations turning $f_T$ into $-f_T$. Now it suffices to
apply Theorem~1-8 from the same volume and Lemma~3.3 and to recall that
$M_T$ is primitive. Then any automorph or anti-automorph of $f_T$ is
of the form $V=\pm (M_T')^n,\ n\in\BZ$, which completes the proof.\qed
\enddemo
\remark{Remark} On the other hand, to prove Proposition~3.4, we may use
Proposition~A.4 (see Appendix).
\endremark

We are going to prove an analog of Theorem~2.6.

\proclaim{Theorem 3.5} Each minimal arithmetic coding of the automorphism $T$
with a primitive matrix $M_T$ is of the form
$$
\phi^\pm_{k,j}(\e)=\(\sum_{n=-\infty}^\infty \e_n\la^{-n}\)
\pmatrix \pm\xi_j\la^k \\ \pm\eta_j\la^k \endpmatrix \mod \BZ^2\tag3.4
$$
for some $n\in\BZ,\ j\in\{0,1,\dots,s\}$. Here $(\xi_j,\eta_j)$
is the homoclinic point given by the solution of equation~{\rm(3.3)}
$(x^{(j)},y^{(j)})$ as follows:
$$
\xi_j=\frac{y^{(j)}+n^{(j)}\la}{\sqrt D},\ \
\eta_j=\frac{x^{(j)}+k^{(j)}\la}{\sqrt D},
$$
and
$$
\pmatrix n^{(j)} \\ k^{(j)}\endpmatrix=
-\det M_T\cdot M_T\pmatrix y^{(j)} \\ x^{(j)}\endpmatrix.
$$
\endproclaim
\demo{Proof} We use practically the same argument as in the proof of
the second part of Theorem~2.6.
Let $\phi_0$ and $\phi$ be two minimal arithmetic codings
for $T$. Recall that the corresponding factor maps
$\phi_0'$ and $\phi'$ are group homomorphisms of the groups
$\X_r'$ (or $\Y_r'$) and $\BT^2$. Suppose first $\Ker\phi'=\Ker\phi'_0$.
Then $A:=\phi\phi_0^{-1}:\BT^2\to\BT^2$ is well defined, and by
definition, $A$ is an automorphism of $\BT^2$ commuting with $T$.
Again, by the Dirichlet theorem and the primitivity of $T$,
we have $A=\pm T^k$, whence $\phi=\pm T^k\phi_0,\ k\in\BZ$.
Thus, for two MAC's with one and the same kernel, the claim is proved.
Since any minimal arithmetic coding is naturally parametrized by
a solution of equation~(3.3), it suffices to apply Lemma~1.1 and
Proposition~3.4.\qed
\enddemo

If $M_T$ is not primitive, this case can be processed in the same spirit;
the corresponding formula for $\phi_{k,j}$ is similar both to formulas~(2.7)
and (3.4).

Following the framework of the previous section (cf. the second part of
Theorem~2.5), we are going to relate minimal arithmetic codings to
the problem of the semiconjungacy of matrices.

\proclaim{Proposition 3.6} Any matrix $B\in GL(2,\BQ)\cap M_2\BZ$ such that
$$
BM_T=C_{r,\sigma} B,\quad \det B=\pm m(T)
$$
has the form
$$
B=\pm\pmatrix x^{(j)}, \, y^{(j)}\\-\det M_T\cdot (x^{(j)},\ y^{(j)})M_T^{-1}
\endpmatrix\cdot M_T^n,\quad n\in\BZ.
$$
Besides,
$\Ker B=\Ker\(\smallmatrix x^{(j)}, \, y^{(j)}\\ (x^{(j)}, \, y^{(j)})
M_T^{-1}\endsmallmatrix\)$, i.e. there is a finite number of possible
kernels for $B$.
\endproclaim
\demo{Proof} A solution $B$ of the matrix equation $BM_T=C_{r,\sigma}B$
together with the condition $\det B=\pm m(T)$ is in fact a matrix
$B=\(\smallmatrix x&y\\z&t\endsmallmatrix\)$, where $(x,y)$ is a
solution of the equation~(3.3), and
$$
(z,t)=-\det M_T\cdot (x,y) M_T^{-1}
$$
(see Theorem~2.5). Now the claim is a direct consequence of
Proposition~3.4.\qed
\enddemo

Thus, we related the problem of description of the kernels of MAC's
for $T$ to the purely algebraic problem of describing the kernels of
the endomorphisms of $\BT^2$ given by the matrices semiconjugating $M$
and $C_{r,\sigma} $. Furthermore, both problems are reduced to
finding the basis solutions of the equation~(3.3). The following example
shows that the situation with distinct series of solutions
can take place, which leads to different series of kernels.

\example{Example} Let $M=\(\smallmatrix 80&9\\9&1\endsmallmatrix\)$.
Then $f_T (x,y)=9x^2-79xy-9y^2$, and it is a direct inspection that
the equations $f_T (x,y)=\pm k$ have no solutions for $1\le k\le8$.
Thus, the integral minimum of $|f_T |$ equals 9. We consider
the equation $f_T (x,y)=\pm9$ and choose the pairs of solutions:
$(x_1=1,y_1=-9)$ and $(x_2=9,y_2=1)$. Constructing now the matrices
$B_1=\(\smallmatrix x_1&y_1\\z_1&t_1\endsmallmatrix\)$ and
$B_2=\(\smallmatrix x_2&y_2\\z_2&t_2\endsmallmatrix\)$ by formula~(2.5),
we obtain thus two matrices from $GL(2,\BQ)\cap M_2\BZ$ semiconjugating
$M$ and the companion matrix $\(\smallmatrix 81&1\\1&0\endsmallmatrix\)$.
However, the matrix $B_1B_2^{-1}$ is not integral, whence the endomorphisms
given by the matrices $B_1$ and $B_2$ have distinct kernels, so do
the corresponding mappings $\phi_1$ and $\phi_2$. Note also that
the groups $\Ker B_1$ and $\Ker B_2$ being isomorphic as
abstract groups, are not isomorphic with respect to $T$ in the sense
that there is no automorphism commuting with $T$ and turning
$\Ker B_1$ into $\Ker B_2$.

Thus, the kernel
of the minimal arithmetic coding is not an invariant for the integral
conjugacy in $GL(2,\BZ)$, as it does not apply even for a single matrix.
\endexample
The idea of this example is based on the fact that $m(T)$ is not a prime.
It can be shown that for $m(T)$ prime such a situation
cannot take place.

\subhead Remark on the case $r<0$ \endsubhead
Finally, we keep our promise and show how to reduce the case $r<0$
to $r>0$. Briefly, given an automorphism $T$ whose matrix $M_T$
has the negative trace, we consider the automorphism with the matrix $-M_T$,
and make sure that it has the same collection of homoclinic
points and the same series with the terms $\e_n\la^{-n}$ but with
$0<\la<1$ and inverted (in the Markov case) restrictions on the digits.

More precisely, let $r<0$ and let $\X_r^-$ be the stationary Markov compactum
$\{\{\e_n\}_{-\infty}^\infty : 0\le\e_n\le |r|,\ \e_n=|r|
\Rightarrow\e_{n-1}=0,\ n\in\BZ\}$. Then any arithmetic coding of $T$ is
given by the mapping
$$
\psi_\t(\e)=\(\sum_{n=-\infty}^\infty\e_n\la^{-n}\)
\pmatrix \xi \\ \eta \endpmatrix\mod \BZ^2,
$$
which formally coincides with the mapping $\phi_\t$
given by formula~(1.3), but acting from $\X_r^-$ if $\sigma=-1$ and $\Y_r$
otherwise with $\la=\frac{r+\sqrt D}2\in(0,1)$ and $(\xi,\eta)$ being
a homoclinic point for $-T$. By formulas~(1.1) and (1.2), the set of
homoclinic points for $T$ and $-T$ is one and the same.
Thus, all claims of the paper for the case $r<0$ remain valid for $r>0$.

Note also that the composition mapping
$S:\X_r\to\X_r^-$ (resp. $\Y_r\to\Y_r$) specified by the formula
$S=\psi_\t^{-1}\phi_\t$ is well defined, does not depend on $\t$,
and $S(\{\e_n\})=\{\e_{-n}\}$.

\head Appendix. Related algebraic questions \endhead

In this appendix we collect all algebraic and number-theoretic
claims which are closely related to the main theorems of the paper, but
at the same time being practically separate. The authors consider them
as known to the specialists or following from certain known facts.
However, some of them prove to be important, namely, Theorem~A.2 which relates
the algebraic conjugacy of the matrices to the equivalence of the
binary quadratic forms,
Theorem~A.7 which answers the question about the number of orbits
of a matrix covering $\BZ^2$, and finally,
Proposition~A.9 describing the Pisot group for a given quadratic PV unit;
we could not find these claims in the classical sources.

\subhead A.1. Unimodular matrices and quadratic forms \endsubhead
We are going to prove an assertion which relates our
theory to the theory of binary integral quadratic forms.
Recall that two binary integral
quadratic forms $f$ and $f'$ are called {\it equivalent} if
$f'(x,y)=f(\a x+\beta y,\ga x +\de y)$ with
$\(\smallmatrix \a&\beta\\\ga&\de\endsmallmatrix\)\in GL(2,\BZ)$.
If $\(\smallmatrix \a&\beta\\\ga&\de\endsmallmatrix\)\in SL(2,\BZ)$,
then they are called {\it properly equivalent}.
For indefinite quadratic forms the problem of equivalence is rather
difficult (see, e.g., \cite{Ven}); note only that for discriminants
appearing in our kind of problems the number of equivalence classes
is large.

Within the appendix we will denote the quadratic
form associated with a matrix $M\in GL(2,\BZ)$, by $f_M$ instead
of $f_T$, which looks more natural here.
Our goal is to prove a claim that relates the problem of the conjugacy of
matrices $M_1$ and $M_2$ in $GL(2,\BZ)$
to the equivalence of the forms accosiated with them.
Let $\goth F$ denote the set of binary integral quadratic forms
with discriminant $r^2\pm4>0$ for some $r\in\BN$. We begin
with a lemma which studies the mapping $\theta:GL(2,\BZ)\to\goth F$
such that $M\overset\theta\to\mapsto f_M$.

\proclaim{Lemma A.1} For a binary quadratic form $f\in\goth F$, the preimage
$\theta^{-1}(f)$ consits exactly of two matrices.
If we denote one of these matrices by $M$, another is
$-\det M\cdot M^{-1}$.
\endproclaim
\demo{Proof} Let $f(x,y)=\a x^2+\beta xy+\ga y^2$, and the matrix sought
is $\(\smallmatrix a&b\\c&d\endsmallmatrix\)=:M$. Then one needs to solve
the equations $\a=b,\beta=d-a,\ga=-c$ for the variables $a,b,c,d$ together
with the condition $ad-bc=\pm1$. Solving them, we see that if
$M_1=\(\smallmatrix a&b\\c&d\endsmallmatrix\)$ is a solution, then another
solution is $M_2=\(\smallmatrix -d&b\\c&-a\endsmallmatrix\)$.\qed
\enddemo

Thus, the mapping $\theta$ is two-to-one, and it is easy to
distinguish the two preimages of a given form, as they have different traces,
though one and the same determinant.

\proclaim{Theorem A.2}
\roster
\item Let matrices $M_1$ and $M_2$ belonging to
$GL(2,\BZ)$ be algebraically conjugate in $GL(2,\BZ)$, i.e.
$BM_1B^{-1}=M_2$ for some $B\in GL(2,\BZ)$.
Then if $\det B=+1$, then the associated forms are properly equivalent.
More precisely, we have $Bf_{M_1}B'=f_{M_2}$. If $\det B=-1$, then
$Bf_{M_1}B'=-f_{M_2}$.
\item If two binary integral quadratic forms $f_1\in\goth F$ and
$f_2\in\goth F$ with one and the same
discriminant are properly equivalent, then
matrices $M_1\in\theta^{-1}(f_1)$ and $M_2\in\theta^{-1}(f_2)$ with
equal traces are algebraically conjugate. If the equivalence
of forms is not proper, then $M_2\sim \det M_1\cdot M_1^{-1}$.
\endroster
\endproclaim
\demo{Proof} (1) Let $M_1=\(\smallmatrix a&b\\c&d\endsmallmatrix\),\
M_2=\(\smallmatrix a'&b'\\c'&d'\endsmallmatrix\)$ with $a+d=a'+d',\
ad-bc=a'd'-b'c'$. Then it is a direct inspection that the relation
$$
B\pmatrix a&b\\c&d\endpmatrix B^{-1}=\pmatrix a'&b'\\c'&d'\endpmatrix
$$
with $B\in SL(2,\BZ)$ is equivalent to the relation
$$
B\pmatrix b&\frac12(d-a)\\\frac12(d-a)&-c\endpmatrix B'=
\pmatrix b'&\frac12(d'-a')\\\frac12(d'-a')&-c'\endpmatrix.
$$
On the contrary, if $B=\(\smallmatrix x&y\\z&t\endsmallmatrix\)
\in GL(2,\BZ)\setminus SL(2,\BZ)$, then the relations $BM_1B^{-1}=M_2$ and
$Bf_{M_1}B'=-f_{M_2}$ in fact lead to the following
one and the same collection of relations:
$$
\align
&f_M(x,y)=-b',\\
&f_M(z,t)=c',\\
&bxz+(d-a)yz-cyt=a-d'.
\endalign
$$
(2) If $B\in SL(2,\BZ)$, then the claim is already proved in the
previous item. If
$B\in GL(2,\BZ)\setminus SL(2,\BZ)$, then, similarly to the first item,
we make sure that the relations $Bf_{M_1}B'=f_{M_2}$
and $BM_1B^{-1}=\det M_2\cdot M_2^{-1}$ also yield one and the same
collection of relations.\qed
\enddemo

The rest of item~A.1 is devoted to diverse applications of this theorem.
The following corollary is straightforward.

\proclaim{Corollary A.3} {\rm(1)} Let two matrices
$M_1$ and $M_2$ from $GL(2,\BZ)$
have one and the same trace and discriminant. Then they are algebraically
conjugate if and only if \newline
-- either the forms $f_{M_1}$ and $f_{M_2}$ are properly equivalent\newline
-- or $f_{M_1}$ is equivalent to $-f_{M_2}$, and the corresponding change
of variables has determinant $-1$.\newline
{\rm(2)} If $\det M=-1$, then $f_M$ is equivalent to $-f_M$. Hence
in this case $M_1\sim M_2$ if and only if $f_{M_1}$ is equivalent
to $f_{M_2}$.
\endproclaim
\demo{Proof} The first item follows from Theorem~A.2. The second one
is a consequence of Lemma~3.3 (recall that $Mf_MM'=-f_M$, if $\det M=-1$).
\qed
\enddemo
The following example shows that sometimes $f_M$ is not equivalent
to $-f_M$ if $\det M=+1$.

\example{Example} Let $M_1=M=\(\smallmatrix 3&5\\1&2\endsmallmatrix\)$,
and $M_2=M^{-1}=\(\smallmatrix 2&-5\\-1&3\endsmallmatrix\)$.
Then $M_1\not\sim M_2$, because the form $f_M(x,y)=5x^2-xy-y^2$ assumes
the value 1, but does not assume the value $-1$.
Indeed, the equation $5x^2-xy-y^2=-1$ is solvable
($x=0,y=1$), while the equation $5x^2-xy-y^2=1$ has no solutions,
as this equation can be rewritten as $(10x-y)^2=2y^2+20$,
whence $\pm2$ should be a quadratic residue modulo 10, what is wrong.
Therefore, $f_M$ is not equivalent to $-f_M$.
\endexample

The first application of Theorem~A.2 is the link between the Dirichlet
theorem for $GL(2,\BZ)$ and the classical theorem on the general form
of a proper automorph of an indefinite binary quadratic form.
We recall that a change of variables with a matrix $B$
is called a {\it proper automorph}
of a form $f$, if $\det B=+1$, and $B'fB=f$.
A form is called {\it primitive}, if its coefficients are relatively
prime.

\proclaim{Proposition A.4} Any proper automorph $B$ of a
primitive binary quadratic
form $f\in\goth F$ is of the form $B=(M')^n,\ n\in\BZ$ for $\det M=+1$ or
of the form $B=(M')^{2n}$ otherwise, where $M\in\theta^{-1}(f)$. The
only exclusion are the forms equivalent to $f_{C_{3,1}}(x,y)=x^2-3xy+y^2$,
whose proper automorph is always of the form $(M')^{n/2}$ for
some $n\in\BZ$.
\endproclaim
\demo{Proof} Let $B$ be a proper automorph of $f$. Choose one of the
matrices in $\theta^{-1}(f)$ and denote it by $M$ (recall that both have the
same determinant). At the first part of the proof of Theorem~2.6 we
in fact proved that if $f$ is primitive, then $M$ is also primitive,
except the exclusive case $r=1,\sigma=-1$
(this link between the terms ``primitive matrix" and ``primitive
form" for completely different notions partially explains our choice
of terminology for the matrices; see also Remark~3 below).

Then by Theorem~A.2, the relation $Bf_MB'=f_M$
implies $BMB^{-1}=M$, whence by the Dirichlet theorem, $B=\pm M^n$,
and it suffices to use the fact that $\det B=+1$. The exclusive case
is studied in the same way.\qed
\enddemo
\remark{Remark {\rm1}} This claim can be obtained by using standard
number-theoretic arguments as a consequence of the
general theorem on the proper automorphs of an indefinite binary quadratic
form, see \cite{Lev, vol. II, Th.~1-8}. To this end, one needs to
find the minimal positive solution of the Pell equation~(0.1) with $+4$.
This way is more computational, while the goal of our proof was to
establish a link with the classical Dirichlet theorem which is applicable
to {\it a priori} completely different class of objects.
\endremark
\remark{Remark {\rm2}} Note that in the nonexclusive case any transformation
of coordinates being either automorph or anti-automorph, is of the
form $\pm(M')^n,\ n\in\BZ$. This is a consequence of the Dirichlet
theorem in its complete form.
\endremark
\remark{Remark {\rm3}} Generally speaking, it is wrong that
the primitivity of a matrix $M$ implies the primitivity of the
associated form $f_M$. Here is the counterexample:
$M=\(\smallmatrix 7&6\\6&5\endsmallmatrix\)$.
\endremark

Now we return to equations~(3.3) (or (2.4)) in order
to find out if the solvability of one of the equations~(3.3)
implies the solvability of another.

\proclaim{Lemma A.5} If $\det M=-1$, then the solvability of one of
the equations~{\rm(3.3)} (say, with $+m$) implies the solvability of another.
On the contrary, for $\det M=+1$ it is, generally speaking, wrong.
\endproclaim
\demo{Proof} If $\det M=-1$, then the claim follows from the equivalence
of the forms $f_M$ and $-f_M$ (see Lemma~3.3) and the fact that
equivalent forms assume one and the same collection of values.
As a couterexample for $\det M=+1$ we can again consider
the matrix $M=\(\smallmatrix 3&5\\1&2\endsmallmatrix\)$ (see Example
above).\qed
\enddemo

Returning now to the problems of Section~2, we will show that
for ``small" discriminants a primitive matrix is always conjugate
to the corresponding companion matrix.

\proclaim{Proposition A.6} (I) Any matrix $M$ with $D=r^2-4\sigma<20$
is algebraically
conjugate in $GL(2,\BZ)$ to the companion matrix $C_{r,\sigma}$.\newline
(II) Let $20\le D<40$ for the matrix $M$ of
an automorphism $T$. Then either $M$ is algebraically conjugate
to the
corresponding companion matrix $C_{r,\sigma}$ or $M$ is not primitive. More
precisely, there is the following alternative.
\roster
\item If $D=21$ or $D=29$, then $M\sim C_{r,\sigma}$.
\item If $D=20$, then either $M$ is primitive and
$M\sim\(\smallmatrix 4&1\\1&0\endsmallmatrix\)$ or
$M\sim\(\smallmatrix 3&2\\2&1\endsmallmatrix\)=
\(\smallmatrix 1&1\\1&0\endsmallmatrix\)^3$.
\item If $D=32$, then either $M\sim C_{6,1}$ or
$M\sim\(\smallmatrix 2&1\\1&0\endsmallmatrix\)^2$.
\endroster
\endproclaim
\demo{Proof} (I) By \cite{Cas1, Ch.~II, \S4, Theorem~VI}, if
positive integers $a,b$ are such that
for an indefinite binary quadratic form $f(x,y)$ with the discriminant $D$
there are no integral $x,y$ such that $-a<f(x,y)<b$, then
$D\ge 4ab+\max\,(a^2,b^2)$. Inverting this assertion and taking
$a=b=2$, we come to the claim of the proposition,
because the integral minimum
of the form $f_M$ in this case equals 1, which is equivalent to
the desired claim.\newline
(II) The central point here is the following sharp estimate for
the integral minimum of an indefinite binary quadratic form: for
a form $f(x,y)=\a x^2+\beta xy +\ga y^2$ with discriminant $D>0$,
$$
\min\,\{|f(x,y)| : (x,y)\in\BZ^2\setminus(0,0)\}\le \sqrt{\frac{25D}{221}}
$$
unless $f$ is equivalent to one of the forms
$l(x^2-xy-y^2)$ or $l(x^2-2y^2)$ with $l\in\BZ$
(see \cite{Cas1, Chap. II, \S4, Th. 6}).
This estimate applied to $f_M$ for $D\le35$ (which actually means that $D<40$)
yields the minimum $=1$, which is equivalent to the solvability
of equation~(2.4). Considering the possibilities for the exclusions
in the cited  claim, we make sure that they could appear only
for $D=20$ or $D=32$. In the first case this leads to the equivalence
of the forms $f_M(x,y)$ and $2(x^2-xy-y^2)$, whence by Theorem~A.2,
$M\sim\(\smallmatrix 3&2\\2&1\endsmallmatrix\)$,
because $\(\smallmatrix 3&2\\2&1\endsmallmatrix\)\sim
\(\smallmatrix 1&-2\\-2&3\endsmallmatrix\)$.
The case $D=32$ is studied in the same way.\qed
\enddemo

\remark{Remark} Practically the first claim of
Proposition~A.6 means that if $r=1,2,3$
for $\sigma=-1$ and $r=3,4$ for $\sigma=+1$, then
$M$ is algebraically conjugate to the corresponding companion matrix.
\endremark
We are going to give some ``counterexamples" showing that the
constants in Proposition~A.6 are precise.

\example{``Counterexamples"} {\bf 1.} The condition $D<20$ cannot be improved.
Indeed, it suffices to consider the matrix
$\(\smallmatrix 3&2\\2&1\endsmallmatrix\)$ which is obviously
not conjugate to the companion matrix. However, this matrix is
not primitive, namely, the cube of $\(\smallmatrix 1&1\\1&0\endsmallmatrix\)$.
\newline
{\bf 2.} For $D=40$ there
exists a primitive matrix with this discriminant not algebraically conjugate 
to the companion matrix, namely,
$M=\(\smallmatrix 5&3\\2&1\endsmallmatrix\)$. Here equation~(2.4) is
$3x^2-4xy-2y^2=\pm1$ and has no integral solutions, as it can be rewritten
as $(3x-2y)^2-10y^2=\pm3$, whence it would follow that $\pm3$ is
a quadratic residue modulo 10.\newline
{\bf 3.} Although each matrix with $D=5$ is algebraically
conjugate by the above to the companion matrix, this is, generally speaking,
wrong for an arbitrary matrix whose spectrum
is in the ring $\BZ[\la]$ with $\la=\frac{1+\sqrt5}2$. Here is an example.
Consider $M=\(\smallmatrix 27&11\\5&2\endsmallmatrix\)$ whose spectrum
is $\{\la^7,-\la^{-7}\}$. A detailed analysis shows that
the integral minimum of the absolute value of the associated form
$f_M(x,y)=11x^2-25xy-5y^2$ equals 5, whence, equation~(2.4) for this case
has no integral solutions, though $M$ is primitive. Thus, it is impossible
to reformulate Proposition~A.6 in purely ``ring" terms.
\endexample
\remark{Remark} Another approach to the problem of conjugacy
of two matrices in $GL(2,\BZ)$ was proposed in \cite{CamTr}.
It is based on the presentation of the group $PSL(2,\BZ)$ as a
free product of cyclic groups. The authors express their gratitude
to B.~Weiss for indicating this reference.
\endremark

\subhead A.2. The number of orbits of a unimodular matrix
\endsubhead

\proclaim{Theorem A.7} Let
$M=\(\smallmatrix a&b\\c&d\endsmallmatrix\)\in GL(2,\BZ)$.
Let $\Orb_M(x,y):=\{M^n\(\smallmatrix x\\y\endsmallmatrix\),\ n\in\BZ\}$
denote the orbit of $(x,y)\in\BZ^2$. Then the linear span of this orbit
$\langle\Orb_M(x,y)\rangle$ is equal to $\BZ^2$ if and
only if
$$
f_M(y,-x)=\pm1.
$$
More generally, for a given $M$,
$$
\min\,\left\{k : \exists \{(x_j,y_j)\}_1^k\mid \bigcup_{j=1}^k
\langle\Orb_M(x_j,y_j)\rangle
=\BZ^2\right\}\ge\min_{(x,y)\neq(0,0)}|f_M(x,y)|.\tag A.1
$$
\endproclaim
\demo{Proof} Similarly to Proposition~3.6, for any
pair $(x,y)\in\BZ^2$ there exists a matrix
$B=B(x,y)=\(\smallmatrix x&z\\y&t\endsmallmatrix\)\in GL(2,\BQ)\cap M_2\BZ$
such that $BC_{r,\sigma}=MB$, where, as above,
$C_{r,\sigma}=\(\smallmatrix r&1\\-\sigma&0\endsmallmatrix\)$ is
the companion matrix. Namely,
$\(\smallmatrix z\\t\endsmallmatrix\)=M^{-1}\(\smallmatrix x\\y
\endsmallmatrix\)$.
Hence $BC_{r,\sigma}^n=M^nB$, and
$$
BC_{r,\sigma}^n\pmatrix 1\\0\endpmatrix=
M^n\pmatrix x\\y\endpmatrix.\tag A.2
$$
Note that by trivial reasons,
$\left\langle C_{r,\sigma}^n\(\smallmatrix 1\\0\endsmallmatrix\),\ n\in\BZ
\right\rangle=\BZ^2$, as this is equivalent to the fact that the powers of
$\la$ form a basis of the module $\BZ[\la]$.

Thus, by relation~(A.2), we have
$$
\langle\Orb_M(x,y)\rangle=B\BZ^2,
$$
whence $\langle\Orb_M(x,y)\rangle$
coincides with $\BZ^2$ if and only if $\det B=\pm1$,
which is equivalent to $\pm1=xt-yz=\pm\bigl(x(-cx+ay)-y(dx-by)\bigr)=
\pm f_M(y,-x)$.

To prove the second claim of the theorem, we observe that
from formula~(A.2) follows the fact that
$\min\,\{|\det B(x,y)|\mid (x,y)\in\BZ^2\setminus\{(0,0)\}\}=
\min\,\{|f_M(x,y)|\mid (x,y)\in\BZ^2\setminus\{(0,0)\}\}=:m$,
whence one needs at least $m$ orbits to cover $\BZ^2$.
\qed
\enddemo
\remark{Remark} Thus, as before, to enumerate all matrices
$B=B(x,y)$ with the minimal possible determinant in modulus,
we need to find all solutions of the Diophantine equation
$$
f_M(y,-x)=\pm m.\tag A.3
$$
Within one and the same ``series of solutions"
of equation~(A.3) (see Proposition~3.4 for the definitions)
$\Orb_M(x,y)\equiv\text{const}$, whence it is easy to construct
an example with the rigid inequality in formula~(A.1).
It suffices to consider any matrix with $m=2$, for instance,
our ``universal" counterexample
$M=\(\smallmatrix 5&3\\2&1\endsmallmatrix\)$. Here all solutions
of equation~(A.3) form a single orbit, hence
the left minimum in the inequality~(A.1) is greater than or equal
to 3.
\endremark

\proclaim{Corollary A.8} If the linear span for the powers of $M$ of
some vector
equals $\BZ^2$, then $M$ is algebraically conjugate to the companion
matrix. Conversely, if $M\sim C_{r,\sigma}$, then there exists
a vector $(x,y)\in\BZ^2$ such that
$\left\langle M^n\(\smallmatrix x\\y\endsmallmatrix\),\ n\in\BZ
\right\rangle=\BZ^2$.
\endproclaim

\subhead A.3. Application to PV numbers \endsubhead
At the end of the appendix we will relate our results to the classical
algebraic theory of Pisot-Vijayaraghavan (PV) numbers.

\definition{Definition} Let $\theta$ be an algebraic integer $>1$ such that
all its Galois conjugates lie inside the unit disc on the complex plane.
Then $\theta$ is called a Pisot-Vijayaraghavan (PV) number.
\enddefinition

Thus, in our case $\la$ is a quadratic PV unit. We define
$$
\Cal P_\la:=\{\xi\in\BR:\|\xi\la^n\|\to0,\quad n\to\infty\}.
$$
Obviously, $\Cal P_\la$ is a group in addition, and from the definition
of a PV number it is clear that if $\la$ is unitary, then
$\BZ[\la]\subset\Cal P_\la$.
We will call $\Cal P_\la$ the {\it Pisot}
group with a parameter $\la$. The implicit description of the Pisot
group is yielded by the classical Pisot-Vijayaraghavan theorem claiming
that $\xi$ belongs to $\Cal P_\la$ if and only if $\xi\in\BQ(\la)$, and
$\Tr(\xi)\in\BZ,\ \Tr(\la\xi)\in\BZ$, where $\Tr(\xi)=\xi+\ov\xi$ is
the trace of a quadratic irrational (see, e.g., \cite{Cas2, Chap.~VIII}).
It is not hard to obtain
the precise description of the Pisot group from these conditions directly,
however, our methods yield its structure almost immediately and relate
this theory to the theory of hyperbolic systems.

\proclaim{Proposition A.9} Let $\la>1$ be the PV number which
satisfies the equation $\la^2=r\la-\sigma$ with $r\ge1$
for $\sigma=-1$ and $r\ge3$ for $\sigma=+1$, and let $D=r^2-4\sigma$. Then
$$
\Cal P_\la=\frac{\BZ+\la\BZ}{\sqrt D}.
$$
\endproclaim
\demo{Proof} Let for simplicity of notation, $\sigma=-1$, and
$\xi\in\Cal P_\la$. Consider the point
$x=(\{\xi\},\{\la^{-1}\xi\})$. Obviously,
$T_{r,\sigma}^n(x)=(\{\la^n\xi\},\{\la^{n-1}\xi\})$,
whence by definition of the Pisot group,
$T_{r,\sigma}^n(x)\to(0,0),\,n\to\pm\infty$,
and thus, $x$ is a homoclinic point for the automorphism $T_{r,\sigma}$,
and by Lemma~2.1, its first coordinate has the form
$\frac{m+n\la}{\sqrt D}\mod 1$ for some $m,n$ integers.

Conversely, let $\xi=\frac{m+n\la}{\sqrt D}$. Then by Lemma~2.1,
the point $x$ is homoclinic for $T_{r,\sigma}$, whence
$\la^n\xi\to0\mod1,\ n\to\infty$.\qed
\enddemo


\Refs

\ref \key Ad
\by R.~L.~Adler
\paper Symbolic dynamics and Markov partitions
\jour Bull. Amer. Math. Soc.
\vol 35
\yr 1998
\pages 1--56
\endref

\ref \key AdWe
\by R. L. Adler and B. Weiss
\paper Entropy, a complete metric invariant for automorphisms of the torus
\jour Proc. Nat. Acad. Sci. USA
\vol 57
\yr 1967
\pages 1573--1576
\endref

\ref \key Ber
\by A.~Bertrand-Mathis
\paper Developpement en base $\theta$, r\'epartition modulo un
de la suite \linebreak
$(x\theta^n)_{n\ge 0}$; langages cod\'es et $\theta$-shift
\jour Bull. Math. Soc. Fr.
\vol 114
\yr 1986
\pages 271--323
\endref

\ref \key BorSh
\by Z. I. Borevich and I. R. Shafarevich
\book Number Theory
\publ New York, Academic Press
\yr 1986
\endref

\ref \key Bow
\by R. Bowen
\paper Markov partitions and measures for axiom A diffeomorphisms
\jour Trans. Amer. Math. Soc.
\vol 154
\yr 1971
\pages 377--397
\endref

\ref \key CamTr
\by J. T. Campbell and E. C. Trouy
\paper When are two elements of $GL(2,\ssize\BZ\dsize)$ similar?
\jour Linear Algebra and its Appl.
\vol 157
\yr 1991
\pages 175--184
\endref


\ref \key Cas1
\by J.~Cassels
\book An introduction to the Geomtery of Numbers
\publ Springer-Verlag, Berlin-G\"ottin\-gen-Heidelberg
\yr 1959
\endref

\ref \key Cas2
\by J.~Cassels
\book An Introduction in Diophantine Approximation
\publ Cambridge Univ. Press
\yr 1957
\endref

\ref \key Fr
\by Ch. Frougny
\paper Representations of numbers and finite automata
\jour Math. Systems Theory
\yr 1992
\vol 25
\pages 37--60
\endref

\ref \key FrSa
\by Ch. Frougny and J. Sakarovitch
\paper Automatic conversion from Fibonacci to golden mean, and generalization
\jour to appear in Int. J. of Alg. and Comput
\endref

\ref \key FrSo
\by Ch. Frougny and B. Solomyak
\paper Finite beta-expansions
\jour Ergod. Theory Dynam. Systems
\vol 12
\yr 1992
\pages 713--723
\endref

\ref \key Gor
\by M. Gordin
\paper Homoclinic approach to the central limit theorem for
dynamical systems
\jour Contemp. Math.
\vol 149
\yr 1993
\pages 149--162
\endref

\ref \key GuSi
\by B. M. Gurevich and Ya. G. Sinai
\paper Algebraic automorphisms of the torus and Markov chains
\publ in book ``P.~Billingsley, Ergodic theory and information
(Russian translation)", Moscow, Mir, 1969
\pages 205--233
\endref

\ref \key KenVer
\by R. Kenyon and A. Vershik
\paper Arithmetic construction of sofic partitions of hyperbolic toral
automorphisms
\jour Erg. Theory Dynam. Systems
\vol 18
\yr 1998
\pages 357--372
\endref

\ref \key Leb
\by S. Le Borgne
\paper Dynamique symbolique et propri\'et\'es stochastiques des
automorphisms du tore : cas hyperbolique et quasi-hyperbolique
\jour Th\`ese doctorale
\yr 1997
\endref

\ref \key Lev
\by W. J. LeVeque
\book Topics in Number Theory
\yr 1956
\publ Addison-Wesley
\endref

\ref \key  Pa
\by W. Parry
\paper On the $\beta$-expansions
of real numbers
\jour Acta Math. Hungar.
\yr 1960
\vol 11
\pages 401--416
\endref

\ref \key SidVer
\by N.~A.~Sidorov and A.~M.~Vershik
\paper Ergodic properties of Erd\"os measure, the entropy of the
golden\-shift, and related problems
\jour to appear in Monatsh. Math
\endref

\ref \key Sin
\by Ya.~Sinai
\paper Markov partitions and A-diffeomorphisms
\jour Funktsional. Anal. i Prilozhen.
\vol 2
\yr 1968
\pages 64--89
\transl\nofrills English transl.
\jour Funct. Anal. Appl.
\vol 2
\yr 1968
\endref

\ref \key Ven
\by B. A. Venkov
\book Elementary Number Theory
\publ Wolters-Noordhoff
\yr 1970
\endref

\ref \key Ver1
\by A. M. Vershik
\paper Locally transversal symbolic dynamics
\jour Algebra i Analiz
\vol 6
\yr 1994
\issue 3
\lang in Russian
\transl \nofrills English transl.
\jour St.~Petersburg Math. J.
\vol 6
\yr 1995
\pages 529--540
\endref

\ref \key Ver2
\by A. M. Vershik
\paper The fibadic
expansions of real numbers and adic
transformation
\inbook Prep. Report Inst.
Mittag--Leffler \yr 1991/1992
\issue 4
\pages 1--9
\endref

\ref \key Ver3
\by A. M. Vershik
\paper Arithmetic isomorphism of the toral hyperbolic automorphisms
and sofic systems
\jour Funktsional. Anal. i Prilozhen.
\lang in Russian
\vol 26
\pages 22--24
\transl\nofrills English transl.
\jour Funct. Anal. Appl.
\yr 1992
\vol 26
\pages 170--173
\endref

\endRefs
\enddocument